\documentclass{article} 

\usepackage{amsmath, amssymb}

\newcommand{\CM}{\mathbb{C}}

\title{ALTERNATING WALK/ZETA CORRESPONDENCE}

\author{Takashi KOMATSU \\
Department of Mathematics, \\
Hiroshima University, \\
Higashihiroshima, Hiroshima 789-8526, JAPAN \\
Math. Research Institute Calc for Industry, \\
Minami, Hiroshima, 732-0816, JAPAN \\ 
e-mail: ta.komatsu@sunmath-calc.co.jp \\
Norio KONNO \\
Department of Applied Mathematics, Faculty of Engineering, \\ 
Yokohama National University \\
Hodogaya, Yokohama 240-8501, JAPAN \\
e-mail: konno-norio-bt@ynu.ac.jp \\ 
Iwao SATO \\ 
Oyama National College of Technology \\
Oyama, Tochigi 323-0806, JAPAN \\ 
e-mail: isato@oyama-ct.ac.jp }
 
\begin{document}
 \maketitle

\clearpage

\begin{abstract}
We consider the alternating zeta function and the alternating $L$-function of a graph $G$, 
and express them by using the Ihara zeta function of $G$. 
Next, we define a generalized alternating zeta function of a graph, and express 
the generalized alternating zeta function of a vertex-transitive regular graph 
by spectra of the transition probability matrix of the symmetric simple random walk 
on it and its Laplacian. 
Furthermore, we present an integral expression for the limit of the generalized 
alternating zeta functions of a series of vertex-transitive regular graphs. 
As an example, we treat the generalized alternating zeta functions of a finite torus. 
Finally, we treat the relation between the Mahler measure and the alternating zeta function 
of a graph. 
\end{abstract}

\vspace{5mm}

{\bf 2000 Mathematical Subject Classification}: 05C50, 15A15. \\
{\bf Key words and phrases} : alternating walk, zeta function, vertex-transitive graph, torus   

\vspace{5mm}

The contact author for correspondence:

Iwao Sato 

Oyama National College of Technology, 
Oyama, Tochigi 323-0806, JAPAN 

E-mail: isato@oyama-ct.ac.jp

\clearpage

\section{Introduction}

Ihara \cite{Ihara} defined the Ihara zeta functions of graphs, and showed that the 
reciprocals of the Ihara zeta functions of regular graphs are explicit polynomials. 
The Ihara zeta function of a regular graph $G$ associated with a unitary 
representation of the fundamental group of $G$ was developed 
by Sunada \cite{Sunada1, Sunada2}. 
Hashimoto \cite{Hashimoto} treated multivariable zeta functions of bipartite graphs. 
Furthermore, Hashimoto  \cite{Hashimoto} gave a determinant expression for the Ihara zeta function 
of a general graph by using the edge matrix. 
Bass \cite{Bass} generalized Ihara's result on the Ihara zeta function of 
a regular graph to an irregular graph $G$. 
Stark and Terras \cite{ST} gave an elementary proof of Bass' Theorem, and 
discussed three different zeta functions of any graph. 
Various proofs of Bass' Theorem were given by Kotani and Sunada \cite{KS}, 
and Foata and Zeilberger \cite{FZ}. 

The Ihara zeta function has many applications in both pure and applied mathematics, 
including for instance dynamical systems, spectral graph theory and complex network analysis. 
In complex network analysis, Arrigo, Higham and Noferini \cite{AHN} introduced a 
non-backtracking alternating walk in a digraph $D$, and presented the exponential
generating function and the resolvent for counting the total number of non-backtracking 
alternating walks of a given length in $D$.
Komatsu, Konnno and Sato \cite{KKS} defined an alternating zeta function of a digraph $D$, and 
presented a determinant expression for this zeta function of $D$. 

Recently, there were exciting developments between quantum walk \cite{Ambainis2003, Kempe2003, 
Kendon2007, Konno2008b, VA} 
on a graph and the Ihara zeta function of a graph. 
We investigated a new class of zeta functions for many kinds of walks including the quantum walk (QW) and 
the random walk (RW) on a graph by a series of ''Zeta Correspondence" of our previous work \cite{K1, K2, K3, K4, K5, K6, K7}. 
In Walk/Zeta Correspondence \cite{K2}, 
a walk-type zeta function was defined without using of the determinant expressions of zeta function 
of a graph $G$, and various properties of walk-type zeta functions of RW, correlated random walk (CRW) 
and QW on $G$ were studied. 
Also, their limit formulas by using integral expressions were presented. 

In this paper, we treat the alternating zeta function of a graph. 
Furthermore, we define a generalized alternating zeta function of a graph, and study various properties for 
the generalized alternating zeta function of a vertex-transitive regular graph. 

This paper is organized as follows:
In Section 2, we give a short review for the Ihara zeta function of a graph. 
In Section 3, we state a review for alternating walk of a digraph. 
In Section 4, we treat the alternating zeta function of a digraph. 
In Section 5, we consider the alternating zeta function of a graph $G$, and express it 
by using the Ihara zeta function of $G$. 
In Section 6, we define a generalized alternating zeta function of a graph, and express 
the generalized alternating zeta function of a vertex-transitive regular graph 
by spectra of the transition probability matrix of the symmetric simple random walk 
on it and its Laplacian. 
In Section 7, we present an integral expression for the limit of the generalized 
alternating zeta functions of a series of vertex-transitive regular graphs. 
In Section 8, we treat the generalized alternating zeta functions of a finite torus.
In Section 9, we deal with the relation between the Mahler measure and the alternating zeta function 
of a graph $G$. 
In Section 10, we consider the alternating zeta function of a regular covering of $G$. 
In Section 11, we express the alternating $L$-function of $G$ by using its Ihara zeta functions. 
In Section 12, we state an example. 

\section{Preliminaries} 

Graphs and digraphs treated here are finite.
Let $G$ be a connected graph and $D(G)= \{ (u,v),(v,u) \mid uv \in E(G) \} $ 
the arc set of the symmetric digraph 
corresponding to $G$. 
For $e=(u,v) \in D(G)$, set $u=o(e)$ and $v=t(e)$. 
Furthermore, let $e^{-1}=(v,u)$ be the {\em inverse} of $e=(u,v)$. 

A {\em path $P$ of length $n$} in $G$ is a sequence 
$P=(e_1, \cdots ,e_n )$ of $n$ arcs such that $e_i \in D(G)$,
$t( e_i )=o( e_{i+1} )(1 \leq i \leq n-1)$. 
If $e_i =( v_{i-1} , v_i )$ for $i=1, \cdots , n$, then we write 
$P=(v_0, v_1, \cdots ,v_{n-1}, v_n )$. 
Set $ \mid P \mid =n$, $o(P)=o( e_1 )$ and $t(P)=t( e_n )$. 
Also, $P$ is called an {\em $(o(P),t(P))$-path}. 
We say that a path $P=( e_1 , \cdots , e_n )$ has a {\em backtracking} 
if $ e^{-1}_{i+1} =e_i $ for some $i(1 \leq i \leq n-1)$. 
A $(v, w)$-path is called a {\em $v$-cycle} 
(or {\em $v$-closed path}) if $v=w$. 
The {\em inverse cycle} of a cycle 
$C=( e_1, \cdots ,e_n )$ is the cycle 
$C^{-1} =( e^{-1}_n, \cdots ,e^{-1}_1 )$.

We introduce an equivalence relation between cycles. 
Two cycles $C_1 =(e_1, \cdots ,e_m )$ and 
$C_2 =(f_1, \cdots ,f_m )$ are called {\em equivalent} if 
$f_j =e_{j+k} $ for all $j$. 
The inverse cycle of $C$ is not equivalent to $C$ if $\mid C \mid \geq 3$. 
Let $[C]$ be the equivalence class which contains a cycle $C$. 
Let $B^r$ be the cycle obtained by going $r$ times around a cycle $B$. 
Such a cycle is called a {\em multiple} of $B$. 
A cycle $C$ is {\em reduced} if 
both $C$ and $C^2 $ have no backtracking. 
Furthermore, a cycle $C$ is {\em prime} if it is not a multiple of 
a strictly smaller cycle. 
Note that each equivalence class of prime, reduced cycles of a graph $G$ 
corresponds to a unique conjugacy class of 
the fundamental group $ \pi {}_1 (G,v)$ of $G$ at a vertex $v$ of $G$. 

The {\em Ihara(-Selberg) zeta function} of $G$ is defined by  
\[
{\bf Z} (G,t)= \prod_{[C]} (1- t^{ \mid C \mid } )^{-1} , 
\]
where $[C]$ runs over all equivalence classes of prime, reduced cycles 
of $G$. 

Let $G$ be a connected graph with $n$ vertices and $m$ edges. 
Then two $2m \times 2m$ matrices 
${\bf B} = {\bf B} (G)=( {\bf B}_{ef} )_{e,f \in D(G)} $ and 
${\bf J}_0 ={\bf J}_0 (G) =( {\bf J}_{ef} )_{e,f \in D(G)} $ 
are defined as follows: 
\[
{\bf B}_{ef} =\left\{
\begin{array}{ll}
1 & \mbox{if $t(e)=o(f)$, } \\
0 & \mbox{otherwise, }
\end{array}
\right.
\   
{\bf J}_{ef} =\left\{
\begin{array}{ll}
1 & \mbox{if $f= e^{-1} $, } \\
0 & \mbox{otherwise.}
\end{array}
\right.
\]
Note that 
\[
{\bf J}_0 = {\bf B} \circ \ {}^t {\bf B} , 
\]
where the {\em Schur/Hadamard product} ${\bf A} \circ {\bf B} $ of two matrices 
${\bf A} $ and ${\bf B} $ is defined by 
\[
( {\bf A} \circ {\bf B} )_{ij} = {\bf A}_{ij} \cdot {\bf B}_{ij} . 
\]

\newtheorem{theorem}{Theorem}
\begin{theorem}[Ihara; Hashimoto; Bass] 
Let $G$ be a connected graph with $n$ vertices and $m$ edges. 
Then the reciprocal of the Ihara zeta function of $G$ is given by 
\[
{\bf Z} (G,t) {}^{-1} = \exp \left(- \sum^{\infty}_{k=1} \frac{N_k }{k} t^k \right)  
= \det ( {\bf I}_{2m} - t ( {\bf B} - {\bf J}_0 )) 
\]
\[
=(1- t^2 )^{m-n} \det( {\bf I}_n  -t {\bf A} (G)+ t^2 ( {\bf D}_G - {\bf I}_n  )) , 
\]
where ${\bf D}_G =( d_{ij} )$ is the diagonal matrix with 
$d_{ii} = \deg {}_G \  v_i \  (V(G)= \{ v_1 , \cdots , v_n \} )$. 
Furthermore, $N_k$ is the number of reduced cycles of length $k$ in $G$ for $k \in \mathbb{N}$. 
\end{theorem}

Let $G=(V(G),E(G))$ be a connected graph and $ x_0 \in V(G)$ 
a fixed vertex.  
Then the {\em generalized Ihara zeta function} $\zeta {}_G (u)$ of $G$ is defined by 
\[
\zeta {}_G (t)= \exp \left( \sum^{\infty}_{m=1} \frac{N^0_m }{m} t^m \right) , 
\]
where $N^0_m $ is the number of reduced $x_0$-cycles of length $m$ in $G$. 
A graph $G$ is called {\em vertex-transitive} if there exists an automorphism $ \phi $ of the automorphism group 
${\rm Aut} (G)$ of $G$ such that $ \phi (u)=v$ for each $u,v \in V(G)$. 
Note that, for a finite vertex-transitive graph, the classical Ihara zeta function is just the above
Ihara zeta function raised to the power equaling the number $n$ of vertices: 
\[
\zeta (G,t)= \zeta (t)={\bf Z} (G,t)^{1/n} . 
\]
Furthermore, the {\em Laplacian} of $G$ is given by 
\[
\Delta = \Delta (G) = {\bf D} - {\bf A} (G). 
\]

A formula for the generalized Ihara zeta function of a vertex-transitive graph is given as follows(see \cite{CJK}):

\begin{theorem}[Chinta, Jorgenson and Karlsson] 
Let $G$ be a  vertex-transitive $(q+1)$-regular graph with spectral measure d $\mu {}_{\Delta }$ for the Laplacian $\Delta $. 
Then 
\[
\zeta (G,t)^{-1} =(1-u^2 )^{(q-1)/2} \exp \left( \int \log (1-(q+1- \lambda )t+q t^2 ) d \mu {}_{\Delta } ( \lambda ) \right) . 
\]
\end{theorem}

\section{An alternating walk of a digraph}   

Let $D$ be a connected digraph with $n$ vertices and $m$ arcs. 
Furthermore, let $A(D)$ be the set of arcs of $D$. 
An {\em alternating walk $P=(v_0, v_1, \cdots ,v_{r-1}, v_r )=[e_1, \cdots ,e_r ]$ of length $r$} in $D$ 
is a sequence of $r$ arcs $e_1, \cdots ,e_r$ such that 
\begin{enumerate} 
\item 
$e_1 =(v_0, v_1), e_2 =(v_2, v_1), \cdots , e_{r-1} =(v_{r-2}, v_{r-1} ), e_r =( v_r , v_{r-1} )$ if $r$ is even, 
\item 
$e_1 =(v_0, v_1), e_2 =(v_2, v_1), \cdots , e_{r-1} =(v_{r-1}, v_{r-2} ), e_r =( v_{r-1} , v_{r} )$ if $r$ is odd, 
\item 
$e_1 =(v_1, v_0), e_2 =(v_1, v_2), \cdots , e_{r-1} =(v_{r-1}, v_{r-2} ), e_r =( v_{r-1} , v_{r} )$ if $r$ is even, 
\item 
$e_1 =(v_1, v_0), e_2 =(v_1, v_2), \cdots , e_{r-1} =(v_{r-2}, v_{r-1} ), e_r =( v_r , v_{r-1} )$ if $r$ is odd.  
\end{enumerate} 
An alternating walk of types 1, 2 (types 3, 4) are called an {\em alternating walk starting with 
an out-edge (an in-edge)}. 
Set $ \mid P \mid =r$, $o(P)=v_0$ and $t(P)= v_r $. 
Also, $P$ is called an {\em $(o(P),t(P))$-alternating walk}. 
An alternating walk $P=[e_1, \cdots ,e_r ]$ of $D$ has {\em backtracking} if $e_{i+1} =e_i $ 
for some $i=1, \ldots, r-1$. 
An alternating walk $P$ is called a {\em non-backtracking alternating walk (NBTAW)} if $P$ has 
no backtracking. 
Otherwise, $P$ is called a {\em backtracking alternating walk (BTAW)}. 

Let ${\bf A}= {\bf A} (D)$ be the adjacency matrix of $D$. 
Furthermore, let an $ n \times n$ matrix $p_k ( {\bf A} )$ denote the matrix whose $(u,v)$-entry 
is the number of NBTAWs from a vertex $u$ to a vertex $v$ of length $k$ that start with an out-edge. 
Moreover, let an $ n \times n$ matrix $q_k ( {\bf A} )$ denote the matrix whose $(u,v)$-entry 
is the number of NBTAWs from a vertex $u$ to a vertex $v$ of length $k$ that start with an in-edge. 
Note that 
\[
p_k ( {\bf A}^T )=q_k ( {\bf A} ) \ for \ all \ k . 
\]

Next, let 
\[ 
{\cal A} = {\cal A} (D)=
\left[ 
\begin{array}{cc}
{\bf 0}_{n} & {\bf A} \\ 
{\bf A}^T & {\bf 0}_n 
\end{array} 
\right] 
, 
\] 
and let 
\[ 
{\bf \Delta} = {\bf \Delta} (D)=
\left[ 
\begin{array}{cc}
{\bf D}_{1} & {\bf 0}_n \\ 
{\bf 0}_n & {\bf D}_2 
\end{array} 
\right] 
, 
\] 
where 
\[ 
{\bf D}_1 = {\bf D}_1 (D)= {\rm diag} ({\bf A} {\bf A}^T) , \ \ {\bf D}_2 ={\bf D}_2 (D)= {\rm diag} ({\bf A}^T {\bf A}) . 
\]
Here, for a square matrix ${\bf F} $, ${\rm diag} ({\bf F} )$ is a diagonal matrix as follows: 
\[
( {\rm diag} ({\bf F} ))_{ij} = 
\left\{
\begin{array}{ll}
{\bf F}_{ii} & \mbox{if $i=j$, } \\
0 & \mbox{otherwise. }
\end{array}
\right. 
\] 
Note that $({\bf D}_1 )_{uu} = {\rm outdeg} \ u$, $({\bf D}_2 )_{uu} = {\rm indeg} \ u$, 
where ${\rm outdeg} \ u$ and ${\rm indeg} \ u$ are the number of arcs having $u$ as an origin and a terminus, respectively.  
Furthermore, we introduce two $2n \times 2n$ matrices $r_{2k} ( {\cal A} )$ and $r_{2k+1} ( {\cal A} )$ are given 
as follows: 
\[
r_{2k} ( {\cal A} )=  
\left[ 
\begin{array}{cc}
p_{2k} ( {\bf A} ) &  {\bf 0}_n \\ 
{\bf 0}_n & q_{2k} ( {\bf A} ) 
\end{array} 
\right] 
, \ \ 
r_{2k+1} ( {\cal A} )=  
\left[ 
\begin{array}{cc}
{\bf 0}_n & p_{2k+1} ( {\bf A} ) \\ 
q_{2k+1} ( {\bf A} ) & {\bf 0}_n 
\end{array} 
\right] 
. 
\]

Arrigo, Higham and Noferini \cite{AHN} gave a two-term recurrence for the matrices $p_k ( {\bf A} )$ and $q_k ( {\bf A} )$, 
and then showed the following result.

\begin{theorem}[Arrigo, Higham and Noferini] 
Let $f(x)= \sum^{\infty}_{k=0} c_k x^k$ be a power series of $x$. 
Furthermore, set $f_h (x)= \sum^{\infty}_{k=0} c_{k+h} x^k$. 
Then 
\[
\sum^{\infty}_{k=0} c_k r_k ( {\cal A} )= 
\left[ 
\begin{array}{cc}
{\bf I}_{2n} & {\bf 0} 
\end{array} 
\right] 
(f_0 ( {\bf Y} )-f_2 ( {\bf Y} ))   
\left[ 
\begin{array}{c}
{\bf I}_{2n} \\
{\bf 0} 
\end{array} 
\right] 
, 
\]  
where 
\[
{\bf Y} =  
\left[ 
\begin{array}{cc}
{\cal A} & {\bf I} - {\bf \Delta} \\ 
{\bf I} & {\bf 0} 
\end{array} 
\right] 
. 
\]
\end{theorem}

Moreover, Arrigo, Higham and Noferini \cite{AHN} gave a formula with respect to $r_{k} ( {\cal A} )$.

\begin{theorem}[Arrigo, Higham and Noferini] 
\[
\left( \sum^{\infty}_{k=0} t^k r_k ( {\cal A} ) \right) ( {\bf I} -t {\cal A} + t^2 ( {\bf \Delta } - {\bf I} ))=
(1- t^2 ) {\bf I} . 
\]
\end{theorem}

Arrigo, Higham and Noferini  \cite{AHN} obtained a resolvent formula with respect to $r_{k} ( {\cal A} )$ 
by Theorem 3. 

Furthermore, we use the Weinstein-Aronszajn identity (see \cite{F}).

\begin{theorem}[the Weinstein-Aronszajn identity]
If ${\bf A}$ and ${\bf B}$ are an $r \times s$ and a $s \times r$ 
matrix, respectively, then we have 
\[
\det ( {\bf I}_{r} - {\bf A} {\bf B} )= 
\det ( {\bf I}_s - {\bf B} {\bf A} ) . 
\] 
\end{theorem}

\section{An alternating zeta function of a digraph}

We introduce an alternating cycle of a digraph $D$, and define a new zeta function 
with respect to alternating cycles of $D$. 

Let $D$ be a connected graph with $n$ vertices $v_1 , \cdots , v_n $ 
and $m$ arcs. 
A $(v, w)$-alternating walk $P=( v= v_0 , v_1 , \ldots , v_r =w)=[ e_1 , \ldots , e_r]$ 
is called a {\em $v$-alternating cycle} (or {\em $v$-closed alternating walk}) if $v=w$ 
and $o( e_1 )=o(e_r )=v$ or $t(e_1 )=t(e_r )=v$. 
Note that the length of each alternating cycle is even. 
The {\em inverse cycle} of an alternating cycle 
$C=[ e_1, \cdots ,e_n ]$ is the alternating cycle 
$C^{-1} =[ e^{-1}_n, \cdots ,e^{-1}_1 ]$.

We introduce an equivalence relation between alternating cycles. 
Two alternating cycles $C_1 =[e_1, \cdots ,e_m ]$ and 
$C_2 =[f_1, \cdots ,f_m ]$ are called {\em equivalent} if 
$f_j =e_{j+k} $ for all $j$. 
The inverse cycle of $C$ is not equivalent to $C$ if $\mid C \mid \geq 3$. 
Let $[C]$ be the equivalence class which contains a cycle $C$. 
An alternating cycle $P=(v_0, v_1, \cdots ,v_{r-1}, v_r )=[e_1, \cdots ,e_r ]$ has a {\em tail} 
if $e_{r} =e_1 $. 
An alternating cycle $C$ is {\em reduced} if $C$ has neither a backtracking nor a tail. 
Let $B^r$ be the alternating cycle obtained by going $r$ times around an alternating cycle $B$. 
Such an alternating cycle is called a {\em multiple} of $B$. 
Furthermore, an alternating cycle $C$ is {\em prime} if it is not a multiple of 
a strictly smaller alternating cycle. 
Then the {\em alternating zeta function} of $D$ is defined by  
\[
{\bf Z}_a (D,t)= \prod_{[C]} (1- t^{ \mid C \mid } )^{-1} , 
\]
where $[C]$ runs over all equivalence classes of prime, reduced alternating cycles 
of $D$.

Now, we consider the following two $m \times m$ matrices ${\bf B}_i={\bf B}_i (D)=( b^{(i)}_{ef} )_{e,f \in A(D)} \ (i=1,2)$: 
\[
b^{(1)}_{ef} =\left\{
\begin{array}{ll}
1 & \mbox{if $o(e)=o(f)$, } \\
0 & \mbox{otherwise, }
\end{array}
\right.
\ 
b^{(2)}_{ef} =\left\{
\begin{array}{ll}
1 & \mbox{if $t(e)=t(f)$, } \\
0 & \mbox{otherwise. }
\end{array}
\right.
\] 
Furthermore, let 
\[
{\cal B} = {\cal B} (D)=   
\left[ 
\begin{array}{cc}
{\bf 0}_m & {\bf B}_1 \\ 
{\bf B}_2 & {\bf 0}_m 
\end{array} 
\right] 
, 
{\bf J} = {\bf J} (D)= 
\left[ 
\begin{array}{cc}
{\bf 0}_m & {\bf I}_m \\ 
{\bf I}_m & {\bf 0}_m 
\end{array} 
\right]  
. 
\] 
Moreover, let $N_k$ be the number of reduced alternating cycles of length $k$ in $D$ 
for $k \in \mathbb{N}$.  
Then the exponential generating function and the determinant expressions for the alternating zeta function 
of a digraph are given as follows.

\begin{theorem}[Komatsu, Konno and Sato]  
Let $D$ be a connected digraph with $n$ vertices and $m$ arcs. 
Then the alternating zeta function of $D$ is given by 
\[
{\bf Z}_a (D,t)= \exp \left( \sum^{\infty}_{k=1} \frac{N_k }{k} t^k \right)  
= \det ( {\bf I}_{2m} -t( {\cal B} - {\bf J} ))^{-1}   
\]
\[
=(1- t^2 )^{-(m-2n)} \det ( {\bf I}_{2n} -t {\cal A} + t^2 ( \Delta - {\bf I}_{2n} ))^{-1} . 
\]
\end{theorem}

The second formula and the third formula are called the {\em Hashimoto expression} and 
the {\em Ihara expression} of ${\bf Z}_a (D,t)$, respectively.

\section{An alternating zeta function of a graph}

We consider the alternating zeta function of the symmetric digraph of a graph. 

Let $G$ be a connected graph with $n$ vertices and $m$ edges, and $D_G $ its symmetric digraph.  
Then we write the alternating function of $D_G$ as follows:
\[
{\bf Z}_a (G,t)= {\bf Z}_a (D_G ,t) . 
\]
We call ${\bf Z}_a (G,t)$ the {\em alternating zeta function} of $G$.  

For the alternating zeta function of a graph, the following result follows.

\begin{theorem}
Let $G$ be a connected graph. 
Then the alternating zeta function of $G$ is given by 
\[
{\bf Z}_a (G,t)= {\bf Z} (G,t) {\bf Z} (G,-t) . 
\]
\end{theorem}

{\em Proof }.  We give two ways of proofs by the Ihara expression and the Hashimoto expression of 
the alternating zeta function. 

I. The proof by the Ihara expression: 

Let $G$ be a connected graph with $n$ vertices and $m$ edges.  
Then we have 
\[
{\bf A} = {\bf A}^T = {\bf A} (G) , \ {\bf D}_{1} = {\bf D}_{2} = {\bf D}_G . 
\]
Thus, we get  
\[ 
\Delta = {\bf I}_2 \otimes {\bf D}_G , \ i.e., \ 
\Delta - {\bf I}_{2n} = {\bf I}_2 \otimes ( {\bf D}_G - {\bf I}_n ) .  
\] 
Set ${\bf Q} = {\bf D}_G - {\bf I}_n $. 
By Theorem 6, we have 
\[
{\bf Z}_a (G,t)^{-1} =(1- t^2 )^{2m-2n} \det ( {\bf I}_{2n} -t {\cal A} + t^2 {\bf I}_2 \otimes {\bf Q} ) .
\]

But, we have 
\begin{align*}
& \det ({\bf I}_{2n} -t {\cal A} + t^2 {\bf I}_2 \otimes {\bf Q} ) \\ 
&= \det 
\left[ 
\begin{array}{cc}
{\bf I}_n + t^2 {\bf Q} & -t {\bf A} \\ 
-t {\bf A} & {\bf I}_n + t^2 {\bf Q}  
\end{array} 
\right]
\ 
\det  
\left[ 
\begin{array}{cc}
{\bf I}_n & t( {\bf I}_n + t^2 {\bf Q} )^{-1} {\bf A} \\ 
{\bf 0} & {\bf I}_n 
\end{array} 
\right]
\\
&= \det 
\left[ 
\begin{array}{cc}
{\bf I}_n + t^2 {\bf Q} & {\bf 0} \\ 
-t {\bf A} & {\bf I}_n + t^2 {\bf Q} - t^2{\bf A} ( {\bf I}_n + t^2 {\bf Q} )^{-1} {\bf A}  
\end{array} 
\right]
\\
&= \det ({\bf I}_n + t^2 {\bf Q} ) 
\det ( {\bf I}_n + t^2 {\bf Q} - t^2 {\bf A} ( {\bf I}_n + t^2 {\bf Q} )^{-1} {\bf A} ) \\ 
&= \det ({\bf I}_n + t^2 {\bf Q} )^2  
\det ( {\bf I}_n - t^2 ( {\bf I}_n + t^2 {\bf Q} )^{-1} {\bf A} ( {\bf I}_n + t^2 {\bf Q} )^{-1} {\bf A} ) \\ 
&= \det ({\bf I}_n + t^2 {\bf Q} )^2  
\det ( {\bf I}_n -t( {\bf I}_n + t^2 {\bf Q} )^{-1} {\bf A} ) \det ( {\bf I}_n +t( {\bf I}_n + t^2 {\bf Q} )^{-1} {\bf A} ) \\ 
&= \det ({\bf I}_n + t^2 {\bf Q} -t {\bf A} ) \det ({\bf I}_n + t^2 {\bf Q} +t {\bf A} ) .   
\end{align*} 
By Theorem 1, we obtain  
\begin{align*}
& {\bf Z}_a (G,t)^{-1} \\  
&= (1- t^2 )^{m-n} \det ({\bf I}_n  -t {\bf A} + t^2 {\bf Q}) \cdot (1- t^2 )^{m-n} \det ({\bf I}_n  +t {\bf A} + t^2 {\bf Q} ) \\  
&= {\bf Z} (G,t) {\bf Z} (G, -t) .   
\end{align*}

II. The proof by the Hashimoto expression: 

At first, let ${\bf K} =( {\bf K}_{ev} )$ ${}_{e \in D(G); v \in V(G)} $ be 
the $m \times n$ matrix defined as follows: 
\[
{\bf K}_{ev} :=\left\{
\begin{array}{ll}
1 & \mbox{if $o(e)=v$, } \\
0 & \mbox{otherwise. } 
\end{array}
\right.
\]
Furthermore, we define the $m \times n$ matrix 
${\bf L} =( {\bf L}_{ev} )_{e \in D(G); v \in V(G)} $ as follows: 
\[
{\bf L}_{ev} :=\left\{
\begin{array}{ll}
1 & \mbox{if $t(e)=v$, } \\
0 & \mbox{otherwise. } 
\end{array}
\right.
\] 
Then we have 
\[
{\bf K} {\bf K}^T = {\bf B}_1 , \ {\bf L} {\bf L}{}^T = {\bf B}_2 , \ {\bf L} {\bf K}^T = {\bf B} .
\]
Moreover, we have 
\[
{\bf K} = {\bf J}_0 {\bf L} , \ {\bf L} = {\bf J}_0 {\bf K} .  
\]

By Theorem 6, we get  
\begin{align*}
{\bf Z}_a (G,t)^{-1} &= \det ({\bf I}_{4m} -t ( {\cal B} - {\bf J} )) \\
&= \det 
\left[ 
\begin{array}{cc}
{\bf I}_{2m} & -t( {\bf B}_1 - {\bf I}_{2m} ) \\ 
-t( {\bf B}_2 - {\bf I}_{2m} ) & {\bf I}_{2m}   
\end{array} 
\right]
\ 
\det  
\left[ 
\begin{array}{cc}
{\bf I}_{2m} & t( {\bf B}_1 - {\bf I}_{2m} ) \\ 
{\bf 0} & {\bf I}_{2m}  
\end{array} 
\right]
\\
&= \det 
\left[ 
\begin{array}{cc}
{\bf I}_{2m} & {\bf 0} \\ 
-t( {\bf B}_2 - {\bf I}_{2m} ) & {\bf I}_{2m} -t^2 ( {\bf B}_2 - {\bf I}_{2m} )( {\bf B}_1 - {\bf I}_{2m} )  
\end{array} 
\right]
\\
&= \det ({\bf I}_{2m} - t^2 ( {\bf B}_2 - {\bf I}_{2m} )( {\bf B}_1 - {\bf I}_{2m} )) \\ 
&= \det ({\bf I}_{2m} - t^2 ( {\bf L} {\bf L}^T - {\bf I}_{2m} )( {\bf K} {\bf K}^T - {\bf I}_{2m} )) . 
\end{align*}

Since 
\[
{\bf J}^2_0 = {\bf I}_{2m} ,  
\]
we have 
\begin{align*}
& {\bf Z}_a (G,t)^{-1} \\
&= \det ({\bf I}_{2m} - t^2 ( {\bf L} {\bf K}^T {\bf J}_0 - {\bf J}^2_0 )( {\bf K} {\bf K}^T - {\bf I}_{2m} )) \\
&= \det ({\bf I}_{2m} - t^2 ( {\bf L} {\bf K}^T - {\bf J}_0 )( {\bf J}_0 {\bf K} {\bf K}^T - {\bf J}_0 )) \\
&= \det ({\bf I}_{2m} - t^2 ( {\bf L} {\bf K}^T - {\bf J}_0 )( {\bf L} {\bf K}^T - {\bf J}_0 )) \\
&= \det ({\bf I}_{2m} - t^2 ( {\bf B} - {\bf J}_0 )^2 ) \\
&= \det ({\bf I}_{2m} -t( {\bf B} - {\bf J}_0 )) \det ({\bf I}_{2m} +t( {\bf B} - {\bf J}_0 )) \\
&= {\bf Z} (G,t) {\bf Z} (G, -t) . 
\end{align*}
$\Box$

Next, we present another proof of Theorem 7 by the Euler product of the alternating zeta function. 
At first, we consider the relation between prime, reduced cycles of a graph $G$ and prime, reduced alternating cycles 
of its symmetric digraph $D_G $.

\newtheorem{proposition}{Proposition}
\begin{proposition}
Let $G$ be a connected graph. 
Then the following result holds: 
\begin{enumerate} 
\item Each prime, reduced cycle with length $2r$ of a graph $G$ corresponds to exactly two prime, reduced alternating cycles 
with length $2r$ of $D_G $; 
\item Each prime, reduced cycle with length $2r+1$ of a graph $G$ corresponds to exactly one prime, reduced alternating cycle  
with length $2(2r+1)$ of $D_G $. 
\end{enumerate} 
\end{proposition}

{\em Proof }.  1: Let $C=( e_1 , e_2 , \ldots , e_{2r-1} , e_{2r} )$ be a prime, reduced cycle with length $2r$ in $G$. 
Then $\tilde{C}=[ e_1 , e^{-1}_2 , \ldots , e_{2r-1} , e^{-1}_{2r} ]$ and $\overline{C}=[ e^{-1}_1 , e_2 , \ldots , e^{-1}_{2r-1} , e_{2r} ]$ 
are prime, reduced alternating cycles with length $2r$ in $D_G$. 

Similarly, the converse is obtained. 

2: Let $C=( e_1 , e_2 , \ldots , e_{2r} , e_{2r+1} )$ be a prime, reduced cycle with length $2r+1$ in $G$. 
Then $\tilde{C}=[ e_1 , e^{-1}_2 , \ldots , e_{2r+1} , e^{-1}_1 , e_2 , \ldots , e^{-1}_{2r+1} ]$ 
is a prime, reduced alternating cycle with length $2(2r+1)$ in $D_G$. 

Similarly, the converse is obtained.  
$\Box$

By Proposition 1, we have the following result.

\begin{theorem}
Let $G$ be a connected graph with $n$ vertices and $m$ edges. 
Then the alternating zeta function of $G$ is given by 
\[
{\bf Z}_a (G,t)= {\bf Z} (G,t) {\bf Z} (G,-t) . 
\]
\end{theorem}

{\em Proof }.  Let ${\cal P}$ be the set of all prime, reduced cycles of $G$. 
Furthermore, let 
\[
{\cal P}_e = \{ [C] \in {\cal P} \mid \ |C| : even \} , \ 
{\cal P}_o = \{ [C] \in {\cal P} \mid \ |C| : odd \} . 
\]
By Proposition 1, we have 
\begin{align*}
& {\bf Z}_a (G,t)^{-1} \\
&= \prod_{[C] \in {\cal P}_e } (1- t^{|C|} )^2 \cdot \prod_{[C] \in {\cal P}_o } (1- t^{2|C|} ) \\ 
&= \prod_{[C] \in {\cal P}_e } (1- t^{|C|} )^2 \cdot \prod_{[C] \in {\cal P}_o } (1- t^{|C|} ) 
\cdot \prod_{[C] \in {\cal P}_o } (1+ t^{|C|} ) \\ 
&= \{ \prod_{[C] \in {\cal P}_e } (1- t^{|C|} ) \prod_{[C] \in {\cal P}_o } (1- t^{|C|} ) \} 
\cdot \{ \prod_{[C] \in {\cal P}_e } (1- t^{|C|} ) \prod_{[C] \in {\cal P}_o } (1+ t^{|C|} ) \} \\
&= {\bf Z} (G,t) {\bf Z} (G, -t) . 
\end{align*}
$\Box$

\section{Alternating Walk/Zeta Correspondence} 

Let $G$ be a connected graph with $n$ vertices and $m$ edges.  
Then the $n \times n$ matrix ${\bf P} = {\bf P} (G)=( P_{uv} )_{u,v \in V(G)}$ is given as follows: 
\[
P_{uv} =\left\{
\begin{array}{ll}
1/( \deg {}_G \ u)  & \mbox{if $(u,v) \in D(G)$, } \\
0 & \mbox{otherwise.}
\end{array}
\right.
\] 
Note that the matrix ${\bf P} (G)$ is the transition probability matrix of the simple RW on $G$. 
By Theorem 7, the alternating zeta function of a regular graph is expressed by the matrix ${\bf P} (G)$ and the Laplacian $\Delta $ of $G$ 
as follows.

\begin{proposition}
Let $G$ be a connected $(q+1)$-regular graph with $n$ vertices and $m$ edges.  
Then 
\begin{align*}
{\bf Z}_a (G,t)^{-1} 
&=(1- t^2)^{2(m-n)} \det ((1+q t^2 ) {\bf I}_n -(q+1)t {\bf P} (G)) \det ((1+q t^2 ) {\bf I}_n +(q+1)t {\bf P} (G)) 
\\
&=(1- t^2)^{2(m-n)}   
\det \left( \{ 1-(q+1)t+q t^2 \} {\bf I}_n +t \Delta \right) \det \left( \{ 1-(q+1)t+q t^2 \} {\bf I}_n -t \Delta \right).  
\end{align*}
\end{proposition}

{\em Proof }.  For the Ihara zeta function of $G$, we have 
\begin{align*}
{\bf Z} (G,t)^{-1} 
&=(1- t^2)^{m-n} \det ((1+q t^2 ) {\bf I}_n -(q+1)t {\bf P} (G)) 
\\
&=(1- t^2)^{m-n}   
\det \left( \{ 1-(q+1)t+q t^2 \} {\bf I}_n +t \Delta \right) .  
\end{align*}
By Theorem 7, the result follows.  
$\Box$

Next, we propose a new zeta function of a graph $G$ and a fixed vertex $x_0 \in V(G)$ .  
Let $G$ be a connected graph.  
Then we define the {\em generalized alternating zeta function} $ \zeta {}_a (G, t)$ of $G$ as follows:    
\[
\zeta {}_a (G, t)= \zeta {}_a (t)= \exp \left( \sum^{\infty}_{k=1} \frac{N^0_k }{k} t^k \right) , 
\]  
where $N^0_k $ is the number of reduced $x_0$-alternating cycles of length $k$ in $G$. 
Note that, for a finite vertex-transitive graph, the generalized alternating zeta function is just the alternating 
zeta function raised to the power equaling the number $n$ of vertices: 
\[
\zeta {}_a (G,t)= {\bf Z}_a (G,t)^{1/n} . 
\]

By Theorem 7, we obtain the following result.

\begin{proposition}
Let $G$ be a connected vertex-transitive $(q+1)$-regular graph with $n$ vertices and $m$ edges.  
Then 
\[
\zeta {}_a (G, t)= \zeta (G,t) \zeta (G, -t) . 
\]
\end{proposition}

{\em Proof }.  By the definition of the generalized alternating zeta function, we have 
\[
\zeta {}_a (G,t)= {\bf Z}_a (G,t)^{1/n} . 
\]
By Theorem 7 and the definition of the generalized Ihara zeta function, we have 
\[
\zeta {}_a (G,t)= {\bf Z} (G,t)^{1/n} {\bf Z} (G,-t)^{1/n} = \zeta (G,t) \zeta (G,-t). 
\]
$\Box$

Now, we present an explicit formula for the generalized alternating zeta function for a vertex-transitive graph. 

Let $G$ be a vertex-transitive $(q+1)$-regular graph with $n$ vertices and $m$ edges. 
Then, since $m=(q+1)n/2$, we have 
\[
\frac{2(m-n)}{n}=q-1 . 
\]

From Proposition 2, we get the following result.

\begin{theorem}[Alternating Walk/Zeta Correspondence] 
Let $G$ be a connected vertex-transitive $(q+1)$-regular graph with $n$ vertices and $m$ edges.   
Then 
\begin{equation} 
\zeta_a (G,t)^{-1} =(1- t^2)^{q-1} 
\exp \left[ \frac{1}{n} \sum_{ \lambda \in {\rm Spec}( {\bf P} )} \log \{ (1+q t^2 )^2 -(q+1)^2 t^2 \lambda {}^2 \} \right] ,
\end{equation} 
\begin{equation}  
\zeta_a (G,t)^{-1} =(1- t^2)^{q-1}  
\exp \left[ \frac{1}{n} \sum_{ \lambda \in {\rm Spec}( \Delta )} \log \left\{ (1-(q+1)t+q t^2 )^2 
- t^2 \lambda {}^2 \right\} \right] . 
\end{equation} 
\end{theorem}

{\bf Proof}.  By Propositions 2, we have 
\begin{align*}
{\zeta}_a (G, t)^{-1} 
&= {\bf Z}_{a} (G, t)^{-1/n} 
\\
&=(1- t^2 )^{2(m-n)/n} \{ \det ((1+q t^2 ) {\bf I}_{n} -(q+1)t {\bf P} ) \det ((1+q t^2 ) {\bf I}_{n} +(q+1)t {\bf P} ) \} {}^{1/n} 
\\
&=(1- t^2 )^{q-1} \left\{ \prod_{ \lambda \in {\rm Spec}( {\bf P} )} ((1+q t^2 )-(q+1)t \lambda )((1+q t^2 )+(q+1)t \lambda ) \right\}^{1/n} 
\\
&=(1- t^2 )^{q-1} \exp \left[ \log \left\{ \prod_{ \lambda \in {\rm Spec}( {\bf P} )} 
((1+q t^2 )^2 -(q+1)^2 t^2 \lambda {}^2 )^{1/n} \right\} \right] 
\\
&=(1- t^2 )^{q-1} \exp \left[ \frac{1}{n} \sum_{ \lambda \in {\rm Spec}( {\bf P} )} \log \{ (1+q t^2 )^2 -(q+1)^2 t^2 \lambda {}^2 \} \right] . 
\end{align*}

Similarly, the second formula follows. 
$\Box$

\section{The generalized alternating zeta functions for the series of regular graphs} 

We present an explicit formula for the generalized alternating zeta functions for the series of regular graphs. 
Let $\{ G_n \}^{\infty}_{n=1} $ be a series of finite vertex-transitive $(q+1)$-regular graphs such that 
\[
\lim {}_{n \rightarrow \infty} |V(G_n )|= \infty . 
\]
Then we have 
\[
\frac{2(|E(G_n )| -|V(G_n )|)}{|V(G_n )|}= \frac{2(q-1)|V(G_n )|}{2|V(G_n )|}=q-1 . 
\]
Set 
\[
\nu_n = |V(G_n)|, \ m_n =|E(G_n )| .
\]

Then the following result holds.

\begin{theorem} 
Let $\{ G_n \}^{\infty}_{n=1} $ be a series of finite vertex-transitive $(q+1)$-regular graphs such that 
\[
\lim {}_{n \rightarrow \infty} |V(G_n )|= \infty . 
\]
Then
\begin{enumerate}   
\item $\lim_{n \rightarrow \infty} \zeta_a (G_n , t)^{-1} =
(1-t^2 )^{q-1} \exp \left[ \int \log \{ (1+q t^2 )^2 -(q+1)^2 t^2 \lambda {}^2 \} d \mu_{P} ( \lambda ) \right] $,   
\item $\lim_{n \rightarrow \infty} \zeta_a (G_n , t)^{-1} =
(1-t^2 )^{q-1} \exp \left[ \int \log \{ (1-(q+1)t+q t^2 )^2 - t^2 \lambda {}^2 \} d \mu_{\Delta} ( \lambda ) \right] $, 
\end{enumerate} 
where $d \mu_P ( \lambda )$ and $d \mu_{\Delta } ( \lambda )$ are the spectral measures for the transition operator ${\bf P}$ 
and the Laplacian $\Delta $. 
\end{theorem}

{\bf Proof}.  By Theorem 9, we have 
\[  
\lim_{n \rightarrow \infty} \zeta_a (G_n , t)^{-1} =
(1-t^2 )^{q-1} \exp \left[ \int \log \{ (1+q t^2 )^2 -(q+1)^2 t^2 \lambda {}^2 \} d \mu {}_P ( \lambda ) \right] .   
\]

Similarly, the second formula follows. 
$\Box$

\section{Torus cases} 

We consider the generalized alternating zeta function of the $d$-dimensional torus $T^d_N \ (d \geq 2)$.  

Let $T^d_N \ (d \geq 2)$ be the {\em $d$-dimensional torus} ({\em graph}) with $N^d$ vertices. 
Its vertices are located in coordinates $i_1 , i_2 , \ldots , i_d $ of a $d$-dimensional Euclidian space $\mathbb{R}^d $, 
where $i_j \in \{ 0,1, \ldots , N-1 \} $ for any $j$ from 1 to $d$. 
A vertex $v$ is adjacent to a vertex $w$ if and only if they have $d-1$ coordinates that are the same, 
and for the remaining coordinate $k$, we have $|i^v_k - i^w_k |=1$, where $i^v_k $ and $i^w_k $ are the $k$-th coordinate 
of $v$ and $w$, respectively. 
Then we have 
\[
|E( T^d_N )|=d N^d ,  
\]
and $T^d_N$ is a vertex-transitive $2d$-regular graph.  

By Proposition 2, we obtain the following result. 
\begin{equation} 
{\bf Z}_a (T^d_N ,t)^{-1} =(1- t^2 )^{2(d-1) N^d} \det ((1+(2d-1) t^2)^2 {\bf I}_{N^d} -4 d^2 t^2 {\bf P} (T^d_N )^2 ) .   
\end{equation} 
Here, it is known that ${\rm Spec}( {\bf P} (T^d_N ))$ is given as follows (see \cite{Spitzer}): 
\[
{\rm Spec}( {\bf P} (T^d_N ))= \left\{ \frac{1}{d} \sum^d_{j=1} \cos \left( \frac{2 \pi k_j }{N} \right) \Bigg|  
k_1 , \ldots , k_d \in \{ 0,1, \ldots , N-1 \} \right\} .
\]
Thus, 
\[  
\displaystyle 
\zeta {}_a (T^d_N ,t)^{-1} =(1- t^2 )^{2(d-1)} \exp \left[ \frac{1}{N^d } 
\sum^d_{j=1} \sum^{N-1}_{ k_j =0} \log \left\{ (1+(2d-1) t^2 )^2 -4 t^2 \left( \sum^d_{j=1} 
\cos \left( \frac{2 \pi k_j }{N} \right) \right)^2 \right\} \right] . 
\]

Therefore, we obtain the following theorem.

\begin{theorem}[Alternating Walk/Zeta Correspondence ($T^d_N$ case)]    
Let $T^d_N \ (d \geq 2)$ be the $d$-dimensional torus with $N^d$ vertices.    
Then 
\[ 
\displaystyle 
\lim_{N \rightarrow \infty} {\zeta}_a (T^d_N ,t)^{-1} =(1- t^2 )^{2(d-1)} \exp \left[ \int^{2 \pi}_{0} \dots \int^{2 \pi}_{0}  
\log \left\{ (1+(2d-1) t^2 )^2 -4 t^2 \left( \sum^d_{j=1} \cos \theta_j \right)^2 \right\} 
\frac{d \theta_1}{2 \pi } \cdots \frac{d \theta_d}{2 \pi } \right] ,  
\] 
where $\int^{2 \pi}_{0} \dots \int^{2 \pi}_{0} $ is the $d$-th multiple integral and 
$ \frac{d \theta_1}{2 \pi } \cdots \frac{d \theta_d}{2 \pi } $ is the uniform measure on $[0, 2 \pi )^d $.  
\end{theorem}

Specially, we consider the case of $d=1$. 
We use the following result (see \cite{K7}).

\newtheorem{lemma}{Lemma} 
\begin{lemma}
Let $r \in \mathbb{R} $ with $|r| \leq 1$. 
Then 
\[
\int_0 ^{2 \pi} \log \left( 1 - r \cdot \sin \theta \right) \frac{d \theta}{2 \pi}
= \int_0 ^{2 \pi} \log \left( 1 - r \cdot \cos \theta \right) \frac{d \theta}{2 \pi}
= \log \left( \frac{1 + \sqrt{1 - r^2}}{2} \right). 
\] 
\end{lemma}

Then the following result follows.

\newtheorem{corollary}{Corollary}
\begin{corollary}     
Let $T^1_N $ be the $1$-dimensional torus with $N$ vertices. 
Then  
\[ 
\displaystyle 
\lim_{N \rightarrow \infty} {\zeta}_a (T^1_N ,t)^{-1} =1. 
\] 
\end{corollary}

{\bf Proof}.  By Theorem 11, we have 
\[
\displaystyle 
\lim_{N \rightarrow \infty} {\zeta}_a (T^1_N ,t)^{-1} = \exp \left[ \int^{2 \pi}_{0} 
\log \left\{ (1+ t^2 )^2 -4 t^2 \cos {}^2 \theta \right\} \frac{d \theta}{2 \pi }  \right] .   
\] 
Thus, 
\begin{align*} 
& \int^{2 \pi}_{0} \log \left\{ (1+ t^2 )^2 -4 t^2 \cos {}^2 \theta \right\} \frac{d \theta}{2 \pi } 
\\
&= \int^{2 \pi}_{0} \left[ \log \left\{ (1+ t^2 )-2t \cos \theta \right\} \frac{d \theta}{2 \pi } 
+ \log \left\{ (1+ t^2 )+2t \cos \theta \right\} \right] \frac{d \theta}{2 \pi } 
\\
&= \int^{2 \pi}_{0} \log \left\{ (1+ t^2 )-2t \cos \theta \right\}  \frac{d \theta}{2 \pi } 
+ \int^{2 \pi}_{0} \log \left\{ (1+ t^2 )+2t \cos \theta \right\} \frac{d \theta}{2 \pi } 
\\
&=2 \log (1+ t^2 )+ \int^{2 \pi}_{0} \log \left\{ 1- \frac{2t}{1+ t^2} \cos \theta \right\} \frac{d \theta}{2 \pi } 
+ \int^{2 \pi}_{0} \log \left\{ 1+ \frac{2t}{1+ t^2} \cos \theta \right\} \frac{d \theta}{2 \pi } . 
\end{align*}

By Lemma 1, we have 
\begin{align*} 
& \int^{2 \pi}_{0} \log \left\{ (1+ t^2 )^2 -4 t^2 \cos {}^2 \theta \right\} \frac{d \theta}{2 \pi } 
\\
&= 2 \log (1+ t^2 )+ \log \left\{ \frac{1+ \sqrt{1-( \frac{2t}{1+ t^2} )^2 }}{2} \right\} 
+ \log \left\{ \frac{1+ \sqrt{1-( \frac{-2t}{1+ t^2} )^2 } }{2} \right\} 
\\
&= 2 \log (1+ t^2 )+ \log \frac{1+t^2 + \sqrt{(1- t^2 )^2 } }{2(1+ t^2 )} 
+ \log \frac{1+t^2 + \sqrt{(1- t^2 )^2 }}{2(1+ t^2 )}  
\\
&= 2 \log \frac{1+t^2 +1- t^2}{2} = 2 \log 1 =0 . 
\end{align*} 
Therefore, the result follows. 
$\Box$

\section{A relation between the Mahler measure and the alternating zeta function} 
  
The {\it logarithmic Mahler measure} $m(f)$ of a nonzero $n$-variable Laurant polynomial $f(X_1, \ldots, X_n)$ is defined by 
\[
m \left( f \right) = \int_{[0,1)^n} \log |f \left( e^{2 \pi i t_1}, \ldots, e^{2 \pi i t_n} \right) | \ d t_1 \cdots  d t_n.
\] 
Note that 
\begin{align*}
m \left( f \right) = \Re \left[ \int_{[0,1)^n} \log \left( f \left( e^{2 \pi i t_1}, \ldots, e^{2 \pi i t_n} \right) \right) \ d t_1 \cdots  d t_n \right],
\end{align*}
where $\Re [z]$ is the real part of $z \in \CM$. 
Sometimes we simply refer to $m(f)$ as the {\it Mahler measure} of $f$. 
Then the above two equalities can be rewritten as 
\begin{align*}
m \left( f \right) 
&= \int_{[0,2 \pi)^n} \log |f \left( e^{i \theta_1}, \ldots, e^{i \theta_n} \right) | \ d \Theta^{(n)}_{unif}
\\
&= \Re \left[ \int_{[0,2 \pi)^n} \log \left( f \left( e^{i \theta_1}, \ldots, e^{i \theta_n} \right) \right) \ d \Theta^{(n)}_{unif} \right] , 
\end{align*} 
where 
\[
d \Theta^{(n)}_{unif} = \frac{d \theta_1}{2 \pi } \cdots \frac{d \theta_n}{2 \pi } . 
\]
This measure was introduced by Mahler \cite{M} in the study of number theory. 
As for Mahler measures, see \cite{GR}, for example. 

Next, we define the {\em logarithmic zeta function} of the generalized alternating zeta function $\zeta {}_a ( T^d_N ,t)$ 
as follows: 
\[
{\cal L} (T^d_{\infty} ,t)= \log \left[ \lim_{N \rightarrow \infty} \left\{ \zeta {}_a (T^d_N ,t)^{-1} \right\} \right] . 
\]

By Theorem 11, we obtain the following result.

\begin{proposition} 
\[  
{\cal L} (T^d_{\infty} ,t)= 2(d-1) \log (1-t^2 ) +  \int^{2 \pi}_{0} \dots \int^{2 \pi}_{0}  
\log \left\{ (1+(2d-1) t^2 )^2 -4 t^2 \left( \sum^d_{j=1} \cos \theta_j \right)^2 \right\} 
\frac{d \theta_1}{2 \pi } \cdots \frac{d \theta_d}{2 \pi } . 
\]
\end{proposition}

Thus,

\begin{theorem}
Let $- \frac{1}{2d-1} <t<0$ and $c=(2d-1)t+t^{-1} $. 
Then 
\begin{align*} 
{\cal L} (T^d_{\infty} ,t) &= 2(d-1) \log (1- t^2 )+2 \log (-t) 
\\
&+ m \left( \sum^d_{j=1} (X_j +X^{-1}_j ) -c \right)+m \left(- \sum^d_{j=1} (X_j +X^{-1}_j ) -c \right) . 
\end{align*} 
\end{theorem}

{\bf Proof}.  By Proposition 6, we have 
\begin{align*} 
{\cal L} (T^d_{\infty} ,t) &= \left[ (d-1) \log (1- t^2 )+ \int^{2 \pi}_{0} \cdots \int^{2 \pi}_{0}  
\log \left\{ 1-2t \left( \sum^d_{j=1} \cos \theta_j \right) +(2d-1) t^2 \right\} 
\frac{d \theta_1}{2 \pi } \cdots \frac{d \theta_d}{2 \pi } \right]  
\\
&+ \left[ (d-1) \log (1- t^2 )+ \int^{2 \pi}_{0} \cdots \int^{2 \pi}_{0}  
\log \left\{ 1+2t \left( \sum^d_{j=1} \cos \theta_j \right) +(2d-1) t^2 \right\} 
\frac{d \theta_1}{2 \pi } \cdots \frac{d \theta_d}{2 \pi } \right] . 
\end{align*}

If $-1<t<0$, then the first term is 
\[
(d-1) \log (1- t^2 )+\log (-t)+  \int^{2 \pi}_{0} \cdots \int^{2 \pi}_{0}  
\log \left\{ \sum^d_{j=1} (e^{i \theta {}_j } + e^{-i \theta {}_j } )-((2d-1)t+ t^{-1} ) \right\} 
\frac{d \theta_1}{2 \pi } \cdots \frac{d \theta_d}{2 \pi } . 
\]
But, $-2d- \{(2d-1)t+ t^{-1} \} >0$ if and only if 
\[
t> - \frac{1}{2d-1} \ or \ t<-1 . 
\]
Thus, if $- \frac{1}{2d-1} <t<0$, then we have 
\[
\sum^d_{j=1} (e^{i \theta {}_j } + e^{-i \theta {}_j } )-((2d-1)t+ t^{-1} ) >0 . 
\]
Therefore, it follows that 
\begin{align*} 
& \int^{2 \pi}_{0} \cdots \int^{2 \pi}_{0} \log \left\{ \sum^d_{j=1} (e^{i \theta {}_j } 
+ e^{-i \theta {}_j } )-((2d-1)t+ t^{-1} ) \right\} \frac{d \theta_1}{2 \pi } \cdots \frac{d \theta_d}{2 \pi } 
\\
&= \Re \left[ \int^{2 \pi}_{0} \cdots \int^{2 \pi}_{0}  
\log \left\{ \sum^d_{j=1} (e^{i \theta {}_j } + e^{-i \theta {}_j } )-((2d-1)t+ t^{-1} ) \right\} 
\frac{d \theta_1}{2 \pi } \cdots \frac{d \theta_d}{2 \pi } \right]  
\\
&= m \left( \sum^d_{j=1} (X_j + X^{-1}_j )-(2(d-1)t +t^{-1} ) \right) . 
\end{align*}

Next, for $- \frac{1}{2d-1} <t<0$, the second term is 
\[
(d-1) \log (1- t^2 )+\log (-t)+  \int^{2 \pi}_{0} \cdots \int^{2 \pi}_{0}  
\log \left\{ - \sum^d_{j=1} (e^{i \theta {}_j } + e^{-i \theta {}_j } )-((2d-1)t+ t^{-1} ) \right\} 
\frac{d \theta_1}{2 \pi } \cdots \frac{d \theta_d}{2 \pi } . 
\] 
Furthermore, if $- \frac{1}{2d-1} <t<0$, then we have 
\[
- \sum^d_{j=1} (e^{i \theta {}_j } + e^{-i \theta {}_j } )-((2d-1)t+ t^{-1} ) >0 . 
\]
Thus,  
\begin{align*} 
& \int^{2 \pi}_{0} \cdots \int^{2 \pi}_{0}  
\log \left\{ - \sum^d_{j=1} (e^{i \theta {}_j } + e^{-i \theta {}_j } )-((2d-1)t+ t^{-1} ) \right\} 
\frac{d \theta_1}{2 \pi } \cdots \frac{d \theta_d}{2 \pi } 
\\
&= \Re \left[ \int^{2 \pi}_{0} \cdots \int^{2 \pi}_{0}  
\log \left\{ - \sum^d_{j=1} (e^{i \theta {}_j } + e^{-i \theta {}_j } )-((2d-1)t+ t^{-1} ) \right\} 
\frac{d \theta_1}{2 \pi } \cdots \frac{d \theta_d}{2 \pi } \right]  
\\
&= m \left( - \sum^d_{j=1} (X_j + X^{-1}_j )-(2(d-1)t +t^{-1} ) \right) . 
\end{align*} 
Therefore, it follows that 
\begin{align*} 
{\cal L} (T^d_{\infty} ,t) &= 2(d-1) \log (1- t^2 )+2 \log (-t) 
\\
&+ m \left( \sum^d_{j=1} (X_j +X^{-1}_j ) -c \right) +m \left( - \sum^d_{j=1} (X_j +X^{-1}_j ) -c \right) . 
\end{align*} 
$\Box$

\section{The alternating zeta function of a regular covering 
of a graph}

Let $G$ be a connected graph, and 
let $N(v)= \{ w \in V(G) \mid (v,w) \in D(G) \} $ denote the 
neighbourhood  of a vertex $v$ in $G$. 
A graph $H$ is called a {\em covering} of $G$ 
with projection $ \pi : H \longrightarrow G $ if there is a surjection
$ \pi : V(H) \longrightarrow V(G)$ such that
$ \pi {\mid}_{N(v')} : N(v') \longrightarrow N(v)$ is a bijection 
for all vertices $v \in V(G)$ and $v' \in {\pi}^{-1} (v) $.
When a finite group $\Pi$ acts on a graph $G$, 
the {\em quotient graph} $G/ \Pi$ is a graph 
whose vertices are the $\Pi$-orbits on $V(G)$, 
with two vertices adjacent in $G/ \Pi$ if and only if some two 
of their representatives are adjacent in $G$.
A covering $ \pi : H \longrightarrow G$ is said to be
{\em regular} if there is a subgroup {\it B} of the 
automorphism group $Aut \  H$ of $H$ acting freely on $H$ such that 
the quotient graph $H/ {\it B} $ is isomorphic to $G$.

Let $G$ be a graph and $ \Gamma $ a finite group.
Then a mapping $ \alpha : D(G) \longrightarrow \Gamma $
is called an {\em ordinary voltage} {\em assignment}
if $ \alpha (v,u)= \alpha (u,v)^{-1} $ for each $(u,v) \in D(G)$.
The pair $(G, \alpha )$ is called an 
{\em ordinary voltage graph}.
The {\em derived graph} $G^{ \alpha } $ of the ordinary
voltage graph $(G, \alpha )$ is defined as follows:
$V(G^{ \alpha } )=V(G) \times \Gamma $ and $((u,h),(v,k)) \in 
D(G^{ \alpha })$ if and only if $(u,v) \in D(G)$ and $k=h \alpha (u,v) $. 
The {\em natural projection} 
$ \pi : G^{ \alpha } \longrightarrow G$ is defined by 
$ \pi (u,h)=u$. 
The graph $G^{ \alpha }$ is called a 
{\em derived graph covering} of $G$ with voltages in 
$ \Gamma $ or a {\em $ \Gamma $-covering} of $G$.
The natural projection $ \pi $ commutes with the right 
multiplication action of the $ \alpha (e), e \in D(G)$ and 
the left action of $ \Gamma $ on the fibers: 
$g(u,h)=(u,gh), g \in \Gamma $, which is free and transitive. 
Thus, the $ \Gamma $-covering $G^{ \alpha }$ is a $ \mid \Gamma \mid $-fold
regular covering of $G$ with covering transformation group $ \Gamma $.
Furthermore, every regular covering of a graph $G$ is a 
$ \Gamma $-covering of $G$ for some group $ \Gamma $ (see \cite{GT}).

Let $G$ be a connected graph, $ \Gamma $ a finite group and 
$ \alpha : D(G) \longrightarrow \Gamma $ an ordinary voltage assignment. 
In the $\Gamma $-covering $G^{ \alpha } $, set $v_g =(v,g)$ and $e_g =(e,g)$, 
where $v \in V(G), e \in D(G), g\in \Gamma $. 
For $e=(u,v) \in D(G)$, the arc $e_g$ emanates from $u_g$ and 
terminates at $v_{g \alpha (e)}$. 
Note that $ e^{-1}_g =(e^{-1} )_{g \alpha (e)}$. 

Let $G$ be a connected graph with $n$ vertices and $m$ edges, 
$ \Gamma $ a finite group and $ \alpha : D(G) \longrightarrow \Gamma $ 
an ordinary voltage assignment. 
For $g \in \Gamma $, let ${\bf A}_g =( a^{(g)}_{uv} )_{u,v \in V(G)} $ be an $n \times n$ matrix as follows: 
\[
 a^{(g)}_{uv} =\left\{
\begin{array}{ll}
1 & \mbox{if $(u,v) \in D(G)$ and $\alpha (u,v)=g$, } \\
0 & \mbox{otherwise.}
\end{array}
\right.
\]
Furthermore, let two $2m \times 2m$ matrices ${\bf B}_{g} ={\bf B}_{g} (G)=( B^{(g)}_{ef} )_{e,f \in D(G)} $ and 
${\bf J}_{g} ={\bf J}_{g} (G)=( J^{(g)}_{ef} )_{e,f \in D(G)} $ are defined as follows: 
\[  
B^{(g)}_{ef} =\left\{
\begin{array}{ll}
1 & \mbox{if $t(e)=o(f)$ and $\alpha (e)=g$, } \\
0 & \mbox{otherwise, }
\end{array}
\right.
\] 
\[  
J^{(g)}_{ef} =\left\{
\begin{array}{ll}
1 & \mbox{if $f= e^{-1} $ and $\alpha (e)=g$, } \\
0 & \mbox{otherwise. }
\end{array}
\right.
\]

Let \( {\bf M}_{1} \oplus \cdots \oplus {\bf M}_{s} \) be the 
block diagonal sum of square matrices 
${\bf M}_{1}, \cdots , {\bf M}_{s}$. 
If \( {\bf M}_{1} = {\bf M}_{2} = \cdots = {\bf M}_{s} = {\bf M} \),
then we write 
\( s \circ {\bf M} = {\bf M}_{1} \oplus \cdots \oplus {\bf M}_{s} \).
The {\em Kronecker product} $ {\bf A} \bigotimes {\bf B} $
of matrices {\bf A} and {\bf B} is considered as the matrix 
{\bf A} having the element $a_{ij}$ replaced by the matrix $a_{ij} {\bf B}$.

A decomposition formula for the Ihara zeta function of 
a regular covering of $G$ is given as follows (\cite{ST2, MS}).

\begin{theorem}[Stark and Terras; Mizuno and Sato] 
Let $G$ be a connected graph with $n$ vertices and $m$ edges, $ \Gamma $ a finite group 
and $ \alpha : D(G) \longrightarrow \Gamma  $ an ordinary voltage assignment. 
Set $ \mid \Gamma \mid =p$. 
Furthermore, let $ {\rho}_{1} =1, {\rho}_{2} , \cdots , {\rho}_{k} $
be the irreducible representations of $ \Gamma $, and 
$d_i$ the degree of $ {\rho}_{i} $ for each $i$, where 
$d_1=1$.
Suppose that the $ \Gamma $-covering $G^{ \alpha } $ of $G$ is connected. 
Then the Ihara zeta function of $G {}^{ \alpha } $ is  
\[
{\bf Z} ( G^{\alpha } , t)^{-1} = {\bf Z} (G,t )^{-1} 
\prod^{k}_{i=2} \det \left( {\bf I}_{2m d_i } -t \left( \sum_{g \in \Gamma } \rho {}_i (g) \otimes {\bf B}_g 
- \sum_{g \in \Gamma } \rho {}_i (g) \otimes {\bf J}_g \right) \right)^{d_i} 
\]
\[
= {\bf Z} (G,t )^{-1} \prod^{k}_{i=2} \left\{ (1- t^2 )^{(m-n) d_i } 
\det \left( {\bf I}_{n d_i } -t \sum_{g \in \Gamma } \rho {}_i (g) \otimes {\bf A}_g 
+ t^2 {\bf I}_{d_i } \otimes {\bf Q} \right) \right\}^{d_i} .   
\]
\end{theorem}

We can generalize the notion of a $\Gamma$-covering of a graph to a simple 
digraph. 
Let $D$ be a connected digraph and $ \Gamma $ a finite group.
Then a mapping $ \alpha : A(D) \longrightarrow \Gamma $
is called a {\em pseudo ordinary voltage} {\em assignment}
if $ \alpha (v,u)= \alpha (u,v)^{-1} $ for each $(u,v) \in A(D)$ 
such that $(v,u) \in A(D)$.
The pair $(D, \alpha )$ is called an 
{\em ordinary voltage digraph}.
The {\em derived digraph} $D^{ \alpha } $ of the ordinary
voltage digraph $(D, \alpha )$ is defined as follows:
$V(D^{ \alpha } )=V(D) \times \Gamma $ and $((u,h),(v,k)) \in 
A(D^{ \alpha })$ if and only if $(u,v) \in A(D)$ and $k=h \alpha (u,v) $. 
The digraph $D^{ \alpha }$ is called a {\em $ \Gamma $-covering} of $D$.
Note that a $\Gamma$-covering of the symmetric digraph corresponding 
to a graph $G$ is a $\Gamma$-covering of $G$(c.f.,  \cite{GT}).

Let $D$ be a connected digraph, $ \Gamma $ a finite group and 
$ \alpha : A(D) \longrightarrow \Gamma $ a pseudo ordinary voltage assignment. 
In the $\Gamma $-covering $D^{ \alpha } $, set $v_g =(v,g)$ and $e_g =(e,g)$, 
where $v \in V(D), e \in A(D), g\in \Gamma $. 
For $e=(u,v) \in A(D)$, the arc $e_g$ emanates from $u_g$ and 
terminates at $v_{g \alpha (e)}$. 
Note that $ e^{-1}_g =(e^{-1} )_{g \alpha (e)}$.

Let $D$ be a connected digraph with $n$ vertices and $m$ arcs, 
$ \Gamma $ a finite group and $ \alpha : A(D) \longrightarrow \Gamma $ 
a pseudo ordinary voltage assignment. 
For $g \in \Gamma $, let ${\bf A}_g =( a^{(g)}_{uv} )_{u,v \in V(D)} $ be an $n \times n$ matrix as follows: 
\[
 a^{(g)}_{uv} =\left\{
\begin{array}{ll}
1 & \mbox{if $(u,v) \in A(D)$ and $\alpha (u,v)=g$, } \\
0 & \mbox{otherwise.}
\end{array}
\right.
\]
Furthermore, let an $m \times m$ matrix ${\bf B}^{\alpha }_{g} ={\bf B}^{\alpha }_{g} (D)=( B^{(g)}_{ef} )_{e,f \in A(D)} $ is defined as follows: 
\[  
B^{(g)}_{ef} =\left\{
\begin{array}{ll}
1 & \mbox{if $t(e)=t(f)$ and $\alpha (e) \alpha (f)^{-1} =g$, } \\
0 & \mbox{otherwise. }
\end{array}
\right.
\]

A decomposition formula for the alternating zeta function of 
a group covering of $D$ is given as follows (\cite{KKS}).

\begin{theorem}[Komatsu, Konno and Sato] 
Let $D$ be a connected digraph with $n$ vertices and $m$ arcs, $ \Gamma $ a finite group 
and $ \alpha : A(D) \longrightarrow \Gamma  $ a pseudo ordinary voltage assignment. 
Set $ \mid \Gamma \mid =p$. 
Furthermore, let $ {\rho}_{1} =1, {\rho}_{2} , \cdots , {\rho}_{k} $
be the irreducible representations of $ \Gamma $, and 
$d_i$ the degree of $ {\rho}_{i} $ for each $i$, where 
$d_1=1$.
Suppose that the $ \Gamma $-covering $D^{ \alpha } $ of $D$ is connected. 
Then the alternating zeta function of $D {}^{ \alpha } $ is  
\[
{\bf Z}_a ( D^{\alpha } , t)^{-1} = {\bf Z}_a (D,t )^{-1} 
\prod^{k}_{i=2} \det ( {\bf I}_{2m d_i } -t( {\cal B}_{\rho {}_i } - {\bf I}_{d_i} \bigotimes {\bf J} ))^{d_i} 
\]
\[
= {\bf Z}_a (D,t )^{-1} \prod^{k}_{i=2} \{ (1- t^2 )^{m d_i -2n d_i } 
\det ( {\bf I}_{2n d_i } -t {\cal A}_{\rho {}_i } + t^2 (( {\bf \Delta } )_{d_i} - {\bf I}_{2n d_i } ) \}^{d_i} ,  
\]
where 
\[
{\cal B}_{\rho {}_i } = 
\left[ 
\begin{array}{cc}
{\bf 0} & {\bf I}_{d_i } \bigotimes {\bf B}_1 \\ 
\sum_{h \in \Gamma } {\rho}_i (h) \bigotimes {\bf B}^{\alpha }_h & {\bf 0}  
\end{array} 
\right] 
, 
{\cal A}_{\rho {}_i } = 
\left[ 
\begin{array}{cc}
{\bf 0} & \sum_{h \in \Gamma } {\rho}_i (h) \bigotimes {\bf A}_h \\ 
\sum_{h \in \Gamma } {\rho}_i (h^{-1} ) \bigotimes {\bf A}^t_h & {\bf 0}
\end{array} 
\right] 
\]
and 
\[
( {\bf \Delta } )_{d_i} = 
\left[ 
\begin{array}{cc}
{\bf I}_{d_i } \bigotimes {\bf D}_{1} & {\bf 0} \\ 
{\bf 0} & {\bf I}_{d_i } \bigotimes {\bf D}_{2} 
\end{array} 
\right] 
. 
\]
\end{theorem}

Now, the alternating zeta function of the $\Gamma $-covering $G^{\alpha } $ of a graph $G$ is given 
as follows:

\begin{theorem} 
Let $G$ be a connected graph with $n$ vertices and $m$ edges, $ \Gamma $ a finite group 
and $ \alpha : D(G) \longrightarrow \Gamma $ an ordinary voltage assignment. 
Set $ \mid \Gamma \mid =p$. 
Furthermore, let $ {\rho}_{1} =1, {\rho}_{2} , \cdots , {\rho}_{k} $
be the irreducible representations of $ \Gamma $, and 
$d_i$ the degree of $ {\rho}_{i} $ for each $i$, where 
$d_1=1$.
Suppose that the $ \Gamma $-covering $G^{ \alpha } $ of $D$ is connected. 
Then the alternating zeta function of $G {}^{ \alpha } $ is  
\[
{\bf Z}_a ( G^{\alpha } , t)^{-1} = {\bf Z} (G,t )^{-1} {\bf Z} (G,-t)^{-1} 
\prod^{k}_{i=2} \det \left( {\bf I}_{2m d_i } 
-t^2 \left( \sum_{g \in \Gamma } \rho {}_i (g) \otimes ( {\bf B}_g - {\bf J}_g ) \right)^2 \right)^{d_i} 
\]
\[
=  {\bf Z} (G,t )^{-1} {\bf Z} (G,-t)^{-1} \prod^{k}_{i=2} \left\{ (1- t^2 )^{2(m-n) d_i } 
\det \left( ( {\bf I}_{n d_i } + t^2 {\bf I}_{d_i} \otimes {\bf Q} )^2 
-t^2 \left( \sum_{g \in \Gamma } \rho {}_i (g) \otimes {\bf A}_g \right)^2 \right) \right\}^{d_i} .   
\]
\end{theorem}

{\em Proof }. By Theorem 7, we have 
\begin{align*}
& {\bf Z}_a (G^{\alpha } ,t)^{-1} = {\bf Z} (G^{\alpha } ,t )^{-1} {\bf Z} (G^{\alpha } ,-t)^{-1} 
\\ 
&= {\bf Z} (G,t )^{-1} 
\prod^{k}_{i=2} \det ( {\bf I}_{2m d_i } -t( \sum_{g \in \Gamma } \rho {}_i (g) \otimes {\bf B}_g 
- \sum_{g \in \Gamma } \rho {}_i (g) \otimes {\bf J}_g ))^{d_i} 
\\
&\ \ \ \ \ \ \ \ \ \ \times {\bf Z} (G,-t )^{-1} 
\prod^{k}_{i=2} \det ( {\bf I}_{2m d_i } +t( \sum_{g \in \Gamma } \rho {}_i (g) \otimes {\bf B}_g 
- \sum_{g \in \Gamma } \rho {}_i (g) \otimes {\bf J}_g ))^{d_i} 
\\ 
&= {\bf Z} (G,t )^{-1} {\bf Z} (G,-t )^{-1} 
\prod^{k}_{i=2} \det ( {\bf I}_{2m d_i } -t^2 ( \sum_{g \in \Gamma } \rho {}_i (g) \otimes ( {\bf B}_g - {\bf J}_g ))^2 )^{d_i} . 
\end{align*}

Furthermore, we get  
\begin{align*} 
& {\bf Z}_a (G^{\alpha } ,t)^{-1} = {\bf Z} (G^{\alpha } ,t )^{-1} {\bf Z} (G^{\alpha } ,-t)^{-1} 
\\ 
&= {\bf Z} (G,t )^{-1} \{ (1- t^2 )^{(m-n) d_i } 
\det ( {\bf I}_{n d_i } -t \sum_{g \in \Gamma } \rho {}_i (g) \otimes {\bf A}_g + t^2 {\bf I}_{d_i } \otimes {\bf Q} ) \}^{d_i} 
\\
&\ \ \ \times {\bf Z} (G,-t )^{-1} \{ (1- t^2 )^{(m-n) d_i } 
\det ( {\bf I}_{n d_i } +t \sum_{g \in \Gamma } \rho {}_i (g) \otimes {\bf A}_g + t^2 {\bf I}_{d_i } \otimes {\bf Q} ) \}^{d_i} 
\\ 
&= {\bf Z} (G,t )^{-1} {\bf Z} (G,-t )^{-1} \prod^{k}_{i=2} \{ (1-t^2 )^{2(m-n) d_i } 
\det (( {\bf I}_{n d_i } + t^2 {\bf I}_{d_i} \otimes {\bf Q} )^2 -t^2 ( \sum_{g \in \Gamma } \rho {}_i (g) \otimes {\bf A}_g )^2 ) \}^{d_i} .  
\end{align*} 
$\Box$

\section{$L$-functions of digraphs}

Let $G$ be a connected graph with $n$ vertices and $m$ edges, 
$ \Gamma $ a finite group and $ \alpha : D(G) \longrightarrow \Gamma $ 
an ordinary voltage assignment. 
Furthermore, let $ \rho $ be a unitary representation of $ \Gamma $ 
and $d$ its degree.
Note that each representation of a finite group is a unitary representation (see \cite{Serre}).

For a cycle $C=( e_1 , e_2 , \ldots , e_{r} )$ of $G$, let 
\[
\rho ( \alpha (C))= \rho ( \alpha (e_1 )) \cdots \rho ( \alpha ( e_{r} )) . 
\] 
The {\em Ihara $L$-function} of $G$ associated with $ \rho $ and $ \alpha $ 
is defined by 
\[
{\bf Z} (G, \rho , \alpha , t)= \prod_{[C]} \det ( {\bf I}_{d} - \rho ( \alpha (C)) t^{|C|} )^{-1} , 
\]
where $[C]$ runs over all equivalence classes of prime, reduced cycle in $G$ (\cite{Ihara, Hashimoto}). 

Let 
\[
N_{\rho , k} = \sum_{C \in {\cal C}_k } \rho ( \alpha (C)) , 
\]
where ${\cal C}_k $ is the set of reduced cycles of length $k$ in $G$. 

Then the following result holds.

\begin{theorem}[Ihara; Hashimoto; Stark and Terrs; Mizuno and Sato] 
Let $G$ be a connected graph with $n$ vertices and $m$ edges, 
$ \Gamma $ a finite group and $ \alpha : D(G) \longrightarrow \Gamma $ 
an ordinary voltage assignment. 
Furthermore, let $ \rho $ be a unitary representation of $ \Gamma $ 
and $d$ its degree. 
Then the reciprocal of the Ihara $L$-function of $G$ associated with 
$ \rho $ and $ \alpha $ is 
\[
{\bf Z} (G, \rho , \alpha , t)= \exp \left( \sum^{\infty}_{k=1} \frac{N_{\rho ,k}}{k} t^k \right) 
= \det \left( {\bf I}_{2md} -t \sum_{g \in \Gamma } \rho (g) \otimes ({\bf B}_g - {\bf J}_g ) \right)^{-1} 
\]
\[
=(1- t^2 )^{-(m-n)d} \det \left( {\bf I}_{nd} -t \sum_{g \in \Gamma } \rho (g) \otimes {\bf A}_g + t^2 {\bf I}_d \otimes {\bf Q} \right)^{-1} .   
\]
\end{theorem}

Let $D$ be a connected digraph with $n$ vertices and $m$ arcs, 
$ \Gamma $ a finite group and $ \alpha : A(D) \longrightarrow \Gamma $ 
a pseudo ordinary voltage assignment. 
Furthermore, let $ \rho $ be a unitary representation of $ \Gamma $ 
and $d$ its degree. 
 
Let $[C]$ be an equivalence class of prime, reduced alternating cycle in $D$. 
Then we consider only an alternating cycle $C=[ e_1 , e_2 , \ldots , e_{2r-1} , e_{2r} ]$ 
such that 
\[
o(C)=o(e_1 )  
\] 
as a representative cycle of $[C]$. 
Furthermore, let 
\[
\tilde{\rho} ( \alpha (C))= \rho ( \alpha (e_1 )) \rho ( \alpha (e_2 ))^{-1} \cdots 
\rho ( \alpha ( e_{2r-1} )) \rho ( \alpha ( e_{2r} ))^{-1} . 
\] 
The {\em alternating $L$-function} of $D$ associated with $ \rho $ and $ \alpha $ 
is defined by 
\[
{\bf Z}_a (D, \rho , \alpha , t)= \prod_{[C]} \det ( {\bf I}_{d} - \tilde{\rho} ( \alpha (C)) t^{|C|} )^{-1} , 
\]
where $[C]$ runs over all equivalence classes of prime, reduced alternating cycle in $D$(\cite{KKS}). 

A determinant expression for the alternating $L$-function of $D$ associated with 
$ \rho $ and $ \alpha $ is given as follows. 
Let 
\[
N_{\rho , k} = \sum_{C \in {\cal C}_k } \tilde{\rho} ( \alpha (C)) , 
\]
where ${\cal C}_k $ is the set of reduced alternating cycles of length $k$ in $D$.

\begin{theorem}[Komatsu, Konno and Sato]
Let $D$ be a connected digraph with $n$ vertices and $m$ arcs, 
$ \Gamma $ a finite group and $ \alpha : A(D) \longrightarrow \Gamma $ 
a pseudo ordinary voltage assignment. 
Furthermore, let $ \rho $ be a unitary representation of $ \Gamma $ 
and $d$ its degree. 
Then the reciprocal of the alternating $L$-function of $D$ associated with 
$ \rho $ and $ \alpha $ is 
\[
{\bf Z}_a (D, \rho , \alpha , t)= \exp \left( \sum^{\infty}_{k=1} \frac{N_{\rho ,k}}{k} t^k \right) 
= \det ( {\bf I}_{2md} -t ( {\cal B}_{ \rho } - {\bf I}_d \bigotimes {\bf J} ) )^{-1} 
\]
\[
=(1- t^2 )^{-(m-2n)d} \det ( {\bf I}_{2nd} -t {\cal A}_{\rho } + t^2 (( {\bf \Delta } )_{d_i} - {\bf I}_{2n d_i } ))^{-1} .   
\]
\end{theorem}

Let $G$ be a connected graph with $n$ vertices and $m$ edges, 
$ \Gamma $ a finite group and $ \alpha : D(G) \longrightarrow \Gamma $ 
an ordinary voltage assignment. 
Furthermore, let $ \rho $ be a unitary representation of $ \Gamma $ 
and $d$ its degree.
Then we write the alternating $L$-function of $D_G$ associated with 
$ \rho $ and $ \alpha $ as follows:
\[
{\bf Z}_a (G, \rho , \alpha ,t)= {\bf Z}_a (D_G , \rho , \alpha ,t) . 
\]
We call ${\bf Z}_a (G, \rho , \alpha ,t)$ the {\em alternating $L$-function} of $G$.  

We obtain a similar result to Theorem 7 for the alternating $L$-function of a graph.

\begin{theorem}
Let $G$ be a connected graph with $n$ vertices and $m$ edges, 
$ \Gamma $ a finite group and $ \alpha : D(G) \longrightarrow \Gamma $ 
an ordinary voltage assignment. 
Furthermore, let $ \rho $ be a unitary representation of $ \Gamma $ 
and $d$ its degree. 
Then the reciprocal of the alternating $L$-function of $G$ associated with 
$ \rho $ and $ \alpha $ is 
\[
{\bf Z}_a (G, \rho , \alpha , t)= {\bf Z} (G, \rho , \alpha , t) {\bf Z} (G, \rho , \alpha , -t) . 
\]
\end{theorem}

{\bf Proof}.  By Theorem 17, we have 
\begin{align*} 
& {\bf Z}_a (D, \rho , \alpha , t)^{-1} 
\\ 
&= \det ( {\bf I}_{2md} -t ( {\cal B}_{ \rho } - {\bf I}_d \bigotimes {\bf J} ) ) 
\\
&= \det  
\left[ 
\begin{array}{cc}
{\bf I}_{2md} & -t ( {\bf I}_d \otimes {\bf B}_1 - {\bf I}_{2md} ) \\ 
-t ( \sum_{h \in \Gamma } {\rho} (h) \bigotimes {\bf B}^{\alpha }_h - {\bf I}_{2md} ) & {\bf I}_{2md}    
\end{array} 
\right]
\\ 
&= \det  
\left[ 
\begin{array}{cc}
{\bf I}_{2md} & -t (  {\bf B}_1 \otimes {\bf I}_d - {\bf I}_{2md} ) \\ 
-t ( \sum_{h \in \Gamma } {\bf B}^{\alpha }_h \bigotimes {\rho} (h) - {\bf I}_{2md} ) & {\bf I}_{2md}   
\end{array} 
\right]
\\ 
&\ \ \ \ \ \ \times  
\det  
\left[ 
\begin{array}{cc}
{\bf I}_{2md} & t ( {\bf B}_1 \otimes {\bf I}_d - {\bf I}_{2md} ) \\ 
{\bf 0} & {\bf I}_{2md}   
\end{array} 
\right]
\\ 
&= \det  
\left[ 
\begin{array}{cc}
{\bf I}_{2md} & {\bf 0} \\ 
-t ( \sum_{h \in \Gamma } {\bf B}^{\alpha }_h \bigotimes {\rho} (h) - {\bf I}_{2md} ) & 
{\bf I}_{2md} - t^2 ( \sum_{h \in \Gamma } {\bf B}^{\alpha }_h \bigotimes {\rho} (h) - {\bf I}_{2md} )
( {\bf B}_1 \otimes {\bf I}_d - {\bf I}_{2md} )    
\end{array} 
\right]
\\ 
&= \det ( {\bf I}_{2md} -t^2 ( \sum_{h \in \Gamma } {\bf B}^{\alpha }_h \bigotimes {\rho} (h) - {\bf I}_{2md} )
( {\bf B}_1 \otimes {\bf I}_d - {\bf I}_{2md} )) .    
\end{align*}

Now, let ${\bf M} =( {\bf M}_{ev} )$ ${}_{e \in D(G); v \in V(G)} $ be the $2md \times nd$ matrix defined 
as follows: 
\[
{\bf M}_{ev} =\left\{
\begin{array}{ll}
\rho ( \alpha (e)) & \mbox{if $t(e)=v$, } \\
{\bf 0}_d & \mbox{otherwise. } 
\end{array}
\right.
\]
Then we have 
\[
{\bf M} \overline{{\bf M}}^T = \sum_{h \in \Gamma } {\bf B}^{\alpha}_h \bigotimes \rho (h) . 
\]
Furthermore, let ${\bf K} =(K_{ev} )_{e \in A(D); v \in V(D)} $ be the $m \times n$ matrix defined as follows:
\[
K_{ev} = \left\{
\begin{array}{ll}
1 & \mbox{if $o(e)=v$, } \\
0 & \mbox{otherwise. } 
\end{array}
\right.
\]
Then we have 
\[
{\bf K} {\bf K}^t = {\bf B}_1 \ and \ {\bf M} ( {\bf K}^T \otimes {\bf I}_d )
= \sum_{g \in \Gamma } {\bf B}_g  \otimes \rho (g) . 
\]

Next, let 
\[
A(D)= \{ e_1 ,\ldots , e_{m} , e^{-1}_1 , \ldots , e^{-1}_{m} \} .   
\]  
Arrange $D(G)$ as follows:
\[
 e_1 ,\ldots , e_{m} , e^{-1}_1 , \ldots , e^{-1}_{m} . 
\] 
Furthermore, let  
\[
{\bf J}_{\rho } = \sum_{g \in \Gamma } {\bf J}_g \otimes \rho (g) . 
\]
Then we have 
\[
{\bf J}_{\rho } = 
\left[ 
\begin{array}{cc}
{\bf 0} & \oplus^m_{i=1} \rho (\alpha (e_i )) \\ 
\oplus^m_{i=1} \rho (\alpha (e^{-1}_i )) & {\bf 0}  
\end{array} 
\right] 
\]
and 
\[
{\bf J}^2_{\rho } = {\bf I}_{2md} \ \ and \ \ \overline{{\bf J}}^T_{\rho } = {\bf J}_{\rho } . 
\]
Thus, we have 
\[
{\bf J}_{\rho }( {\bf K} \otimes {\bf I}_d )= {\bf M} . 
\]
Therefore, it follows that 
\begin{align*} 
& {\bf Z}_a (D, \rho , \alpha , t)^{-1} 
\\ 
&= \det ( {\bf I}_{2md} -t^2 ( {\bf M} \overline{{\bf M}}^T - {\bf I}_{2md} )
( {\bf K} {\bf K}^T \otimes {\bf I}_d - {\bf I}_{2md} )) 
\\ 
&= \det ( {\bf I}_{2md} -t^2 ( {\bf M} ( {\bf K}^T \otimes {\bf I}_d ) \overline{{\bf J}}^T_{\rho } - {\bf J}_{\rho } \overline{{\bf J}}^T_{\rho })
( {\bf K} {\bf K}^T \otimes {\bf I}\d - {\bf I}_{2md} )) 
\\
&= \det ( {\bf I}_{2md} -t^2 ( {\bf M} ( {\bf K}^T \otimes {\bf I}_d )- {\bf J}_{\rho } )
( {\bf J}_{\rho } ( {\bf K} {\bf K}^T \otimes {\bf I}_d )- {\bf J}_{\rho } )) 
\\ 
&= \det ( {\bf I}_{2md} -t^2 ( {\bf M} ( {\bf K}^T \otimes {\bf I}_d )- {\bf J}_{\rho } )
(  {\bf M} ( {\bf K}^T \otimes {\bf I}_d ) - {\bf J}_{\rho } )) 
\\ 
&= \det ( {\bf I}_{2md} -t ( {\bf M} ( {\bf K}^T \otimes {\bf I}_d )- {\bf J}_{\rho } ))  
\det ( {\bf I}_{2md} +t ( {\bf M} ( {\bf K}^T \otimes {\bf I}_d )- {\bf J}_{\rho } )) 
\\ 
&= \det ( {\bf I}_{2md} -t \sum_{g \in \Gamma } \rho (g) \otimes ( {\bf B}_g - {\bf J}_g )) 
\det ( {\bf I}_{2md} +t \sum_{g \in \Gamma } \rho (g) \otimes ( {\bf B}_g - {\bf J}_g )) 
\\ 
&= {\bf Z} (G, \rho , \alpha , t)^{-1} {\bf Z} (G, \rho , \alpha , -t)^{-1} . 
\end{align*} 
$\Box$

By Theorems 14 and 17, we obtain the following result.

\begin{corollary}
Let $G$ be a connected graph, $ \Gamma $ a finite group and
$ \alpha : D(G) \longrightarrow \Gamma  $ an ordinary voltage assignment. 
Suppose that the $ \Gamma $-covering $G^{ \alpha } $ of $G$ is connected. 
Then 
\[
{\bf Z}_a ( G^{\alpha } , t )= \prod_{ \rho } {\bf Z}_a (G, \rho , \alpha , t)^{ \deg \rho} , 
\]
where $\rho $ runs over all inequivalent irreducible representation of $\Gamma $. 
\end{corollary}

\section{Example}

Let $D$ be the symmetric digraph of the complete graph $G=K_3$ with vertices 1,2,3. 
Then we have 
\[
A(D)= \{ a , b, c,  a^{-1} , b^{-1} , c^{-1} \} , 
\]
where $a=(1,2), \ b=(2,3) , c=(3,1)$. 
The equivalence classes of prime, reduced cycles of $G$ are 
\[
\{ (a,b,c), (b,c,a), (c,a,b) \}, 
\{ (a^{-1} , c^{-1} , b^{-1} ), (c^{-1} , b^{-1} , a^{-1} ), (b^{-1} , a^{-1} , c^{-1} ) \} . 
\]
Furthermore, the equivalence classes of prime, reduced alternating cycles of $G$ are $[C] , \ [C^{-1} ]$, 
where 
\[
C=[ a, b^{-1} , c, a^{-1}, b, c^{-1} ] , \ C^{-1} =[ c, b^{-1} , a, c^{-1}, b, a^{-1} ] . 
\]
By the definition of the alternating function, we have 
\[
{\bf Z} (G,t)^{-1} =(1-t^3 )^2, \ \ {\bf Z}_a (G,t)^{-1} = (1- t^{|C|} ) (1- t^{|C^{-1} |} )=(1- t^6 )^2 . 
\]

But, 
\[
{\bf A}=  
\left[ 
\begin{array}{ccc}
0 & 1 & 1 \\
1 & 0 & 1 \\
1 & 1 & 0 
\end{array} 
\right]
, \ 
{\bf A}^T = {\bf A} 
, \ 
{\bf D}_1 = {\bf D}_2 =2 {\bf I}_3 . 
\]
By Theorem 7, we have 
\begin{align*}
{\bf Z}_a (G,t)^{-1} &= {\bf Z} (G,t)^{-1} \cdot {\bf Z} (G, -t)^{-1} \\
&= (1-t^3 )^2 (1+t^3)^2 =(1-t^6 )^2 . 
\end{align*}

Next, let $ \Gamma = \mathbb{Z}_3 = \{ 1, \tau , \tau {}^2 \} ( \tau {}^3 =1)$ be 
the cyclic group of order 3 and $ \alpha : D(K_3 ) \longrightarrow \mathbb{Z}_3$ 
the ordinary voltage assignment such that 
$ \alpha (1,2)= \tau , \alpha (2,3)=1$ and $ \alpha (3,1)=1$.
Then the $Z_3$-covering $G^{ \alpha } $ is isomorphic to the symmetric digraph of the 
cycle graph $C_9$ of length 9, and so there exist only two equivalence classes of 
prime, reduced alternating cycles of length 18 in $G^{\alpha } $. 
By the definition of the alternating function of a digraph, we have 
\[
{\bf Z}_a (D^{\alpha } ,t)^{-1} =(1- t^{18} )^2 . 
\]

The characters of ${\bf Z}_3$ are given as follows:
\[
\chi {}_i ( \tau {}^j )=( \eta {}^i )^j , \  0 \leq i,j \leq 2, \  
\eta = \exp \left( \frac{2 \pi \sqrt{-1} }{3} \right) = \frac{-1+ \sqrt{-3} }{2} . 
\]
Then we have 
\[
{\bf A}_1 =  
\left[ 
\begin{array}{ccc}
0 & 0 & 1 \\
0 & 0 & 1\\
1 & 1 & 0 
\end{array} 
\right]
, \ 
{\bf A}_{\tau } =  
\left[ 
\begin{array}{ccc}
0 & 1 & 0 \\
0 & 0 & 0 \\
0 & 0 & 0 
\end{array} 
\right]
, \  
{\bf A}_{\tau {}^2 } =  
\left[ 
\begin{array}{ccc}
0 & 0 & 0 \\
1 & 0 & 0 \\
0 & 0 & 0 
\end{array} 
\right]
. 
\]

By Theorem 16, we have 
\begin{align*} 
{\bf Z} (G, \chi {}_1 , \alpha , t)^{-1} &= (1-t^2 )^{3-3} 
\det ( {\bf I}_3 -t \sum^2_{i=0} \chi {}_1 ( \tau {}_i ) {\bf A} {}_{\tau {}_i } +t^2 ( {\bf D}_G - {\bf I}_3 )) 
\\
&= \det ( {\bf I}_3 -t \sum^2_{i=0} \chi {}_1 ( \tau {}_i ) {\bf A} {}_{\tau {}_i } +t^2 {\bf I}_3 ) .  
\end{align*} 
Thus, we have 
\begin{align*} 
{\bf Z} (G, \chi {}_1 , \alpha , t)^{-1} &= \det ( 
\left[ 
\begin{array}{ccc}
1+ t^2 & - \eta t & -t \\
- \eta^2 t & 1+t^2 & -t \\
-t & -t & 1+t^2 
\end{array} 
\right]  
) 
\\
&= (1+t^2 )^3 -( \eta + \eta {}^2 ) t^3 -3t^2 (1+t^2)=1+ t^3 +t^{6} .  
\end{align*} 
By Theorem 18, we have 
\begin{align*} 
{\bf Z}_a (G, \chi {}_1 , \alpha , t)^{-1} &= {\bf Z} (G, \chi {}_1 , \alpha , t)^{-1} {\bf Z} (G, \chi {}_1 , \alpha , -t)^{-1} 
\\ 
&= (1+t^3 +t^6 )(1- t^3 +t^6 )=1+ t^6 +t^{12} . 
\end{align*}

Similarly, we have 
\[ 
{\bf Z}_a (G, \chi {}_2 , \alpha , t)^{-1} = {\bf Z}_a (G, \chi {}_1 , \alpha , t)^{-1} 
=1+ t^6 +t^{12} . 
\]
By Corollary 2, we have 
\begin{align*} 
{\bf Z}_a ( G^{\alpha } , t)^{-1} &= {\bf Z}_a ( G, t)^{-1} {\bf Z}_a ( G, \chi {}_1 , \alpha , t)^{-1} 
{\bf Z}_a ( G, \chi {}_2 , \alpha , t)^{-1} 
\\
&= (1-t^6 )^2 (1+t^6 +t^{12} )^2 = (1-t^{18} )^2 . 
\end{align*}

{\bf Acknowledgments} 
 
The first author is supported by the JSPS Grant-in-aid for young scientisits No. 22K13959.


\begin{thebibliography}{99}


\bibitem{Ambainis2003} 
A. Ambainis,  
Quantum walks and their algorithmic applications, 
Int. J. Quantum Inf. {\bf 1} (2003), 507--518. 

\bibitem{AAR} 
G. E. Andrew, R. Askey and R. Roy, 
Special Functions, 
Cambridge University Press, VCambridge (1999). 

\bibitem{AHN}
F. Arrigo, D. J. Higham and V. Noferini, 
Non-backtracking alternating walks, 
SIAM Journal on Applied Mathematics 79, (2019), 781-801. 

\bibitem{Bass}
H. Bass,
The Ihara-Selberg zeta function of a tree lattice,
Internat. J. Math. 3 (1992), 717-797.

\bibitem{CJK} 
G. Chinta, J. Jorgenson and A. Karlsson,  
Heat kernels on regular graphs and generalized Ihara zeta function formulas,  
Monatsh. Math. {\bf 178} (2015), 171-190. 

\bibitem{F} 
H. Flanders, 
The elementary divisors of $AB$ and $BA$, 
Proc. Am. Math. Soc. 2 (1951), 871-874. 

\bibitem{FZ}
D. Foata and D. Zeilberger, 
A combinatorial proof of Bass's evaluations of the Ihara-Selberg zeta 
function for graphs, 
Trans. Amer. Math. Soc. 351 (1999), 2257-2274. 

\bibitem{GT}
J. L. Gross and T. W. Tucker,
Topological Graph Theory,
Wiley-Interscience, New York, 1987.

\bibitem{GR} 
A. J. Guttmann and M. D. Rogers, 
Spanning tree generating functions and Mahler measure, 
J. Phys. A: Math. Theor. {\bf 43} (2010), 305205. 

\bibitem{Hashimoto}
K. Hashimoto,
Zeta Functions of Finite Graphs and Representations 
of $p$-Adic Groups,
Adv. Stud. Pure Math. Vol. 15, Academic Press, 
New York, 1989, pp. 211-280. 

\bibitem{Ihara} 
Y. Ihara,
On discrete subgroups of the two by two projective linear group 
over $p$-adic fields,
J. Math. Soc. Japan 18 (1966), 219-235.

\bibitem{Kempe2003} 
J. Kempe,  
Quantum random walks - an introductory overview, 
Contemporary Physics {\bf 44} (2003), 307--327.  

\bibitem{Kendon2007} 
V. Kendon, 
Decoherence in quantum walks - a review, 
Math. Struct. in Comp. Sci. {\bf 17} (2007), 1169--1220. 

\bibitem{KKS} 
T. Komatsu, N. Konno and I. Sato,  
A zeta function with respect to non-backtracking alternating walks for a digraph, 
Linear Algebra and its Applications {\bf 620} (2021), 344-367. 

\bibitem{K1}
T. Komatsu, N. Konno and I. Sato,  
Grover/Zeta Correspondence based on the Konno-Sato theorem, 
Quantum Inf. Process. {\bf 20}, 268 (2021). 

\bibitem{K2}
T. Komatsu, N. Konno and I. Sato,  
Walk/Zeta Correspondence, 
arXiv:2104.10287 (2021).   

\bibitem{K3}
T. Komatsu, N. Konno and I. Sato,  
IPS/Zeta Correspondence, 
Quantum Inf. Comput. {\bf 22} (2022), 251--269. 

\bibitem{K4}
T. Komatsu, N. Konno and I. Sato,  
Vertex-Face/Zeta Correspondence,  
J. Algebraic Combin. (in press), arXiv:2107.03300 (2021).  

\bibitem{K5}
T. Komatsu, N. Konno and I. Sato,  
CTM/Zeta Correspondence,  
Quantum Stud.: Math. Found. {\bf 9} (2022), 165--173. 

\bibitem{K6}
T. Komatsu, N. Konno, I. Sato and S. Tamura,  
A Generalized Grover/Zeta Correspondence,  
arXiv:2201.03973 (2022).  

\bibitem{K7}
T. Komatsu, N. Konno, I. Sato and S. Tamura,  
Mahler/Zeta Correspondence,   
to appear in QIP. 

\bibitem{Konno2008b} 
N. Konno,  
Quantum Walks, 
In: Lecture Notes in Mathematics: Vol.1954, pp.309--452, Springer-Verlag, Heidelberg (2008). 

\bibitem{K7}
N. Konno and S. Tamura, 
Walk/Zeta Correspondence for quantum and correlated random walks,  
Yokohama Math. J. {\bf 67} (2021), 125--152. 

\bibitem{KS}
M. Kotani and T. Sunada,
Zeta functions of finite graphs,
J. Math. Sci. U. Tokyo 7 (2000) 7-25.

\bibitem{M} 
K. Mahler, 
On some inequalities for polynomials in several variables, 
J. Lond. Math. soc. {\bf 37} (1962), 341-3444. 

\bibitem{MS} 
H. Mizuno and I. Sato, 
Zeta functions of graph coverings, 
Journal of Combin. Theory, Ser. B {\bf 80} (2000), 247--257. 

\bibitem{Serre}
J. -P. Serre,
Linear Representations of Finite Group,
Springer-Verlag, New York, 1977.

\bibitem{Spitzer}  
F. Spitzer,  
Principles of Random Walk (2nd edition). 
New York (NY): Springer; 1976. 

\bibitem{ST}
H. M. Stark and A. A. Terras, 
Zeta functions of finite graphs and coverings, 
Adv. Math. 121 (1996), 124-165. 

\bibitem{ST2} 
H. M. Stark and A. A. Terras, 
Zeta functions of finite graphs and coverings, part II, 
Advances in Math. {\bf 154} (2000), 132--195. 

\bibitem{Sunada1}
T. Sunada,
$L$-Functions in Geometry and Some Applications,
in Lecture Notes in Math., Vol. 1201, Springer-Verlag, New York, 1986, 
pp. 266-284. 

\bibitem{Sunada2}
T. Sunada,
Fundamental Groups and Laplacians (in Japanese), 
Kinokuniya, Tokyo, 1988. 

\bibitem{VA} 
S. E. Venegas-Andraca,     
Quantum walks: a comprehensive review,   
Quantum Inf. Process. {\bf 11} (2012), 1015-1106. 


\end{thebibliography}
\end{document}